\documentclass[12pt]{article}
\usepackage{amsthm,amsfonts,amssymb,amscd}

\begin{document}
\begin{center}
\Large \bf Birationally rigid varieties \\
with a pencil of Fano double covers. II
\end{center}
\vspace{0.7cm}

\centerline{\large \bf Aleksandr V. Pukhlikov}
\vspace{0.7cm}

\begin{center}
Max-Planck-Institut f\" ur Mathematik \\
Vivatsgasse 7 \\
53111 Bonn \\
GERMANY \\
e-mail: {\it pukh@mpim-bonn.mpg.de}
\end{center}
\vspace{0.3cm}

\begin{center}
Steklov Institute of Mathematics \\
Gubkina 8 \\
117966 Moscow \\
RUSSIA \\
e-mail: {\it pukh@mi.ras.ru}
\end{center}
\vspace{0.3cm}

\begin{center}
Division of Pure Mathematics \\
Department of Mathematical Sciences \\
M$\&$O Building, Peach Street \\
The University of Liverpool \\
Liverpool L69 7ZL \\
ENGLAND\\
e-mail: {\it pukh@liv.ac.uk}
\end{center}
\vspace{0.7cm}

\centerline{March 1, 2004}\vspace{0.7cm}


\parshape=1
3cm 10cm \noindent {\small \quad \quad \quad
\quad\quad\quad\quad\quad\quad\quad {\bf Abstract}\newline We
continue to study birational geometry of Fano fibrations
$\pi\colon V\to {\mathbb P}^1$ the fibers of which are Fano double
hypersurfaces of index 1. For a majority of families of this type,
which do not satisfy the condition of sufficient twistedness over
the base, we prove birational rigidity (in particular, it means
that there are no other structures of a fibration into rationally
connected varieties) and compute their groups of birational
self-maps. We considerably improve the principal components of the
method of maximal singularities, in the first place, the technique
of counting multiplicities for the fibrations $V/{\mathbb P}^1$
into Fano varieties over the line.} \vspace{0.7cm}

\newpage

CONTENTS \vspace{0.5cm}

\noindent Introduction

0.1. $K^2$-condition and $K$-condition

0.2. The list of varieties under consideration

0.3. Varieties with a pencil of double spaces

0.4. Formulation of the main result

0.5. The structure of the paper

0.6. Historical remarks

0.7. Acknowledgements \vspace{0.3cm}

\noindent 1. The method of maximal singularities

1.1. Maximal singularities of linear systems

1.2. A stronger version of the Noether-Fano inequality

1.3. The self-intersection of the linear system $\Sigma$

1.4. $K^2$-condition and birational rigidity

1.5. The generalized $K^2$-condition \vspace{0.3cm}

\noindent 2. The technique of counting multiplicities

2.1. The notations and the principal claim

2.2. Proof of Proposition 2.1: counting multiplicities

2.3. Proof of Lemma 2.1

2.4. Estimating the multiplicities of a linear system:

the non-singular case

2.5. Estimating the multiplicities of a linear system:

the singular case \vspace{0.3cm}

\noindent 3. Varieties with a pencil of Fano double covers

3.1. Movable systems on the varieties of the type ((0),(2,0))

3.2. Checking the $K$-condition

\parshape=1
3cm 10cm \noindent {\small 3.2.1. Varieties of the type
((0),(1,1)) \newline 3.2.2. Varieties of the type ((2),(1,0))
\newline 3.2.3. Varieties of the type ((2),(1,0)) \newline 3.2.4.
Varieties of the type ((2),(0,0)) \newline 3.2.5. Varieties of the
type ((3),(0,0)) \newline 3.2.6. Varieties of the type
((1,2),(0,0)) \newline 3.2.7. Varieties of the type
((1,1,1),(0,0)}

3.3. Proof of birational rigidity

3.4. Multiplicities of subvarieties of codimension 2

3.5. Estimating the number of lines

3.6. A method of estimating the degree \vspace{0.3cm}

\noindent References

\newpage

\section*{Introduction}

The present paper is a direct follow up of the paper [16]. We
study birational geometry of higher-dimensional algebraic
varieties with a pencil of Fano double covers, now without the
assumption that the standard condition of sufficient twistedness
over the base, that is, the $K^2$-condition [7,10], holds. As we
pointed out in [10], if the deviation from the $K^2$-condition is
not too big, the techniques of the method of maximal singularities
still works and makes it possible to prove birational rigidity.
This is the subject of the present paper: we consider Fano
fibrations $V/{\mathbb P}^1$ which do not satisfy the
$K^2$-condition, but which, however, present a not too strong
deviation from this condition.

If the deviation from the $K^2$-condition oversteps a certain
boundary, then the variety $V$ is no longer birationally rigid.
The non-rigid and ``boundary'' families of fibrations $V/{\mathbb
P}^1$ will be considered in the next paper, the third part of this
research project. For that purpose we will use the improved
technique developed here.

\subsection{$K^2$-condition and $K$-condition}

In this paper, as in the previous paper [16], we deal with Fano
fibrations $\pi\colon V\to{\mathbb P}^1$, satisfying the
conditions
$$
A^1V=\mathop{\rm Pic}V={\mathbb Z}K_V\oplus {\mathbb Z}F,\quad
A^2V={\mathbb Z}K^2_V\oplus{\mathbb Z}H_F,
$$
where $F$ is the class of the projection $\pi$, $H_F=(-K_V\cdot
F)$ is the hyperplane section of the fiber. Set $A^1_{\mathbb
R}V=A^1V\otimes{\mathbb R}$, $A^2_{\mathbb R}V=A^2V\otimes{\mathbb
R}$. Let
$$
A^1_+V\subset A^1_{\mathbb R}V,\quad A^1_{\mathop{\rm
mov}}V\subset A^1_{\mathbb R}V,\quad A^2_+V\subset A^2_{\mathbb
R}V
$$
be the closed cones, generated respectively by effective divisors,
movable divisors and effective cycles of codimension two. Geometry
of the fibration $V/{\mathbb P}^1$ is to a considerable extent
determined by the position of the class $K^2_V$ with respect to
the cone $A^2_+V$, and also by the position of the anticanonical
class $(-K_V)$ with respect to the cone $A^1_{\mathop{\rm mov}}V$.
Obviously, $A^2_{\mathbb R}V={\mathbb R}K^2_V\oplus{\mathbb
R}H_F$, so that $A^2_{\mathbb R}V^*={\mathbb R}\beta\oplus{\mathbb
R}\chi$, where $\beta\colon A^2_{\mathbb R}V\to{\mathbb R}$ is
defined by the condition $\beta(H_F)=1,\beta(K^2_V)=0$ and
similarly $\chi(K^2_V)=1,\chi(H_F)=0$.


It is well known (and easy to prove, see [7,10]) that the
$K^2$-condition
$$
K^2_V\not\in \mathop{\rm Int}A^2_+V
$$
implies the $K$-condition:
$$
-K_V\not\in \mathop{\rm Int}A^1_{\mathop{\rm mov}}V.
$$
(There is a natural {\it self-intersection map}
$$
\mathop{\rm sq}\colon A^1_{\mathop{\rm mov}}V\to A^2_+V,
$$
$$
\mathop{\rm sq}\colon z\mapsto z^2,
$$
and it is easy to check that $\mathop{\rm sq}(\mathop{\rm
Int}A^1_{\mathop{\rm mov}}V)\subset\mathop{\rm Int}A^2_+V$, see
Remark 1.1 in \S 1 below.) The present paper deals with Fano
fibrations $V/{\mathbb P}^1$, each fiber of which is a regular
Fano double hypersurface of index 1, $F_t=\pi^{-1}(t)$,
$F_t\in{\cal F}^{\mathop{\rm reg}}$, see [16], whereas the
$K^2$-condition does not hold, that is, $K^2_V\in \mathop{\rm
Int}A^2_+V$. The deviation from the $K^2$-condition is measured by
the number $a\geq 0$, satisfying the formula
$$
K^2_V-aH_F\in \partial A^2_+V.
$$
In this paper we consider mainly the fibrations $V/{\mathbb P}^1$,
satisfying the $K$-condition. The $K$-condition is also not
absolutely necessary for birational rigidity: if the deviation
from the $K$-condition is not too big, the techniques of the
method of maximal singularities still work and make it possible to
complete the study. The families that do not satisfy the
$K$-condition and also certain other families which we do not
consider here and which were not considered in [16] will be
studied in the next paper, the third part of this research.

\subsection{The list of varieties under consideration}

Recall the construction of varieties with a pencil of Fano double
hypersurfaces, see [9,16] for details. Let ${\cal
E}=\bigoplus{\cal O}_{{\mathbb P}^1}(a_i)$ be a locally free sheaf
of rank $M+2$, where $a_0=0\leq a_1\leq\dots\leq a_{M+1}$, $X
={\mathbb P}({\cal E})$, $\mathop{\rm Pic}X={\mathbb Z}L_X
\bigoplus{\mathbb Z}R$, where $L_X$ is the class of the
tautological sheaf, $R$ is the class of a fiber of the natural
projection $\pi_X\colon X\to{\mathbb P}^1$. Now the variety $V$ is
realized as the double cover $\sigma\colon V\to Q$ of the smooth
hypersurface $Q\subset X$,
$$
Q\sim mL_X+a_QR,\quad a_Q\in{\mathbb Z}_+,
$$
branched over a smooth divisor $W_Q=W\cap Q$, where $W\subset X$,
$$
W\sim 2lL_X+2a_WR,\quad a_W\in{\mathbb Z}_+.
$$
By the symbol $\pi\colon V\to{\mathbb P}^1$ we denote the natural
projection, by the symbol $F_t$ the fiber $\pi^{-1}(t)$, sometimes
omitting $t$. It is easy to see that
$$
K_V=-L_V+(a_X+a_Q+a_W-2)F,
$$
where $L_V=\sigma^*(L_X|_Q)$, $a_X=a_1+\dots+a_{M+1}$. Since the
linear system $|L_X|$ (and therefore $|L_V|$) is free, the
$K^2$-condition follows from the inequality
\begin{equation}
\label{i1} (K^2_V\cdot L^{M-1}_V)=2m(4-a_X-2a_Q-2a_W)+2a_Q\leq 0.
\end{equation}
In this paper we consider families that do not satisfy (\ref{i1}).
Let us give their list. The parameters of these families are
written in the following format
$$
((a_1,\dots,a_{M+1}),(a_Q,a_W)),
$$
and for brevity of notations in the set $(a_1,\dots,a_{M+1})$ we
write only the non-zero entries, if there are any, otherwise we
write $(0)$: thus, for example, (1) means the set $(0,\dots,0,1)$,
$(1,1)$ means the set $(0,\dots,0,1,1)$, and $(0)$ stands for the
set $(0,\dots,0)$, consisting entirely of zeros. Here is the list
of varieties studied in this paper:
$$
\begin{array}{cl}
1.&((0),(2,0))\\
2.&((0),(1,1))\\
3.&((1),(0,1))\\
4.&((2),(1,0))\\
5.&((2),(0,0))\\
8.&((3),(0,0))\\
7.&((1,2),(0,0))\\
8.&((1,1,1),(0,0))\\
\end{array}
$$
Once again we emphasize that the eight families listed above do
not complete the list of varieties that do not satisfy the $K^2$
condition. There are seven more families not satisfying the
$K$-condition, most of which are not birationally rigid, having
many structures of a rationally connected fibration. These
families are not considered in this paper.

\subsection{Varieties with a pencil of double spaces}

Birational rigidity of varieties $V$, fibered into double spaces
of index one, was proved in the very first paper of the author,
devoted to Fano fibrations, in [7], in the assumption of
sufficient twistedness over the base, that is, the
$K^2$-condition. In this paper we also consider the varieties of
this type, satisfying the $K$-condition but not the
$K^2$-condition. Recall their construction, which is a
simplification of the construction for double hypersurfaces given
above (corresponding to the value $m=1$). Let ${\cal
E}=\bigoplus{\cal O}_{{\mathbb P}^1}(a_i)$ be a locally free sheaf
of rank $M+1$, $a_0=0\leq\dots\leq a_M$, $X={\mathbb P}(\cal E)$
its projective bundle, then $V$ is realized as the double cover
$\sigma\colon V\to X$, branched over a smooth hypersurface
$W\subset X$,
$$
W\sim 2ML_X+2a_WR,
$$
where $L_X$ is the class of the tautological sheaf, $R$ is the
class of a fiber of the projection $\pi_X\colon X\to {\mathbb
P}^1$. Obviously,
$$
K_V=-L_V+(a_X+a_W-2)F,
$$
where $L_V=\sigma_*L_X$, $F=\sigma^*R$ is the class of a fiber,
$$
\mathop{\rm Pic}V={\mathbb Z}L_V\bigoplus{\mathbb Z}F,
$$
$a_X=a_1+\dots+a_M,a_W\in{\mathbb Z}_+$. The parameters of a
family are written in the following format:
$$
((a_1,\dots,a_M),a_W),
$$
where as above if $(a_1,\dots,a_M)\neq(0,\dots,0)$, then the zeros
are omitted, and if $(a_1,\dots,a_M)=(0,\dots,0)$, then we write
simply $(0)$. A variety $V$ satisfies the $K^2$-condition, if the
following inequality holds:
\begin{equation}
\label{i2} (K^2_V\cdot L^{M-1}_V)=8-2a_X-4a_W\leq 0.
\end{equation}
In the present paper we consider the families that do not satisfy
(\ref{i2}). Here is their list:
$$
\begin{array}{cl}
1^*.&((1),1)\\
2^*.&((2),0)\\
3^*.&((3),0)\\
4^*.&((1,2),0)\\
5^*.&((1,1,1),0)\\
\end{array}
$$
For all the types listed above the $K$-condition is satisfied (it
will be proved below). There are four families more, for which
both the $K^2$-condition and $K$-condition are not satisfied.
Their birational geometry will be studied in the next paper.

\subsection{Formulation of the main result}

For a divisor $D$ on a rationally connected variety $Y$ we set
$$
c(D,Y)=\mathop{\rm sup}\{\varepsilon\in {\mathbb
Q}_+|D+\varepsilon K_Y\in A^1_+Y\}.
$$
For a movable linear system $\Sigma$ on a uniruled variety $V$
define the {\it virtual threshold of canonical adjunction} by the
formula
$$
c_{virt}(\Sigma)=\mathop{\rm inf}\limits_{V^{\sharp}\to V}
\{c(\Sigma^{\sharp},V^{\sharp})\},
$$
where the infimum is taken over all birational morphisms
$V^{\sharp}\to V$, where $V^{\sharp}$ is a projective model of the
field ${\mathbb C}(V)$, non-singular in codimension 1,
$\Sigma^{\sharp}$ is the strict transform of the system $\Sigma$
on $V^{\sharp}$. The following definition of birational rigidity
is equivalent to the standard one.

{\bf Definition 0.1.} (i) A variety $V$ is said to be {\it
birationally superrigid}, if for any movable linear system
$\Sigma$ on $V$ the following equality holds:
$$
c_{\mathop{\rm virt}}(\Sigma)=c(\Sigma,V).
$$
(ii) A variety $V$ (respectively, a Fano fibration $V/S$) is said
to be {\it birationally rigid}, if for any movable linear system
$\Sigma$ on $V$ there exists a birational self-map
$\chi\in\mathop{\rm Bir}V$ (respectively, a fiber-wise self-map
$\chi\in\mathop{\rm Bir}(V/S)$), satisfying the following
equality:
$$
c_{\mathop{\rm virt}}(\Sigma)=c(\chi_*\Sigma,V).
$$

Let us formulate the main result of the paper.

{\bf Theorem 1.} (i) {\it Regular varieties $V$ with a pencil of
double hypersurfaces of the types $2-8$ in the list above and all
varieties of the types $1^*-5^*$ with a pencil of double spaces
are birationally superrigid and satisfy the $K$-condition:
$-K_V\not\in A^1_{\mathop{\rm mov}}V$. The defined structure of a
rationally connected fibration $V/{\mathbb P}^1$ is the unique
non-trivial structure of a rationally connected fibration on $V$.
The groups of birational and biregular self-maps of these
varieties coincide:}
$$
\mathop{\rm Bir}V=\mathop{\rm Aut}V={\mathbb Z}/2{\mathbb Z}
$$

(ii) {\it The regular varieties $V$ of the type 1 (that is, the
varieties of the family $((0),(2,0))$) are birationally rigid.
They have a non-trivial birational self-map, an involution
$\tau\in\mathop{\rm Bir}V\setminus\mathop{\rm Aut}V$, and do not
satisfy the $K$-condition. On $V$ there are exactly two
non-trivial structures of a rationally connected fibration: the
pencil $|F|$ of fibers of the morphism $\pi$ and its image
$|\tau_*F|$. The group $\mathop{\rm Bir}V$ consists of four
elements:}
$$
\mathop{\rm Bir}V=<\tau>\times\mathop{\rm Aut}V=({\mathbb
Z}/2{\mathbb Z})^{\times 2}.
$$

{\bf Corollary 0.1.} {\it All varieties of the types $1-8$ and
$1^*-5^*$ do not admit a structure of a rationally connected
fibration over a base of dimension 2 or higher, in particular,
they cannot be fibered by a rational map into rational curves or
rational surfaces. Therefore they all are non-rational.}

\subsection{The structure of the paper}

To prove Theorem 1, it is necessary to improve the technique of
the method of maximal singularities. Such an improvement is
important by itself since it extends the domain where the method
works: birational rigidity is proved for far from all the natural
families of Fano fibrations over ${\mathbb P}^1$, even in the
assumption of sufficient twistedness over the base. For this
reason, a bigger part of the paper is of technical character.

In \S 1 we carry out some preparatory work: we prove that the
virtual and actual thresholds of canonical adjunction coincide
provided that the $K^2$-condition holds and the fibers of the Fano
fibration $V/{\mathbb P}^1$  satisfy the standard conditions (h),
(v) and (vs) (see [16]). After that we make the first step in the
direction of improving the technique: we show that replacing the
condition (h) by a stronger one compensates replacing the
$K^2$-condition by a weaker one.

In \S 2 we improve the technique of the method of maximal
singularities on the basis of a radically different idea: to
compare singularities of the horizontal cycle $Z^h$ and of its
restriction onto the fiber $F$, containing the centre of the
maximal singularity. This idea is new. We use it to prove
birational rigidity for the first time. Besides we give in full
detail the technique of counting multiplicities: in all the
previous papers [7-16] we computed the multiplicities of the cycle
$Z$ of intersection of two divisors of the linear system $\Sigma$,
whereas here we intersect a divisor and an irreducible subvariety
of codimension 2. The computations are parallel to the divisorial
case, however there are some differences: for example, only the
blow ups of subvarieties of codimension 3 and higher are taken
into account whereas the blow ups of subvarieties of codimension 3
play the same part as the blow ups of subvarieties of codimension
2 in the divisorial case.

In \S 3 we prove Theorem 1. First we check the $K$-condition for
varieties from the lists above. After that the necessary estimates
are verified for subvarieties $Y\subset F$ of an arbitrary fiber
$F$ (this work was mostly carried out in the previous paper [16]
and in [9]). The remaining part of the section contains an
improvement of one estimate for multiplicities of subvarieties of
codimension two. In order to do that, one needs to get an estimate
for the number of lines (counted with multiplicities) through an
arbitrary point of a fiber. Combining this improved estimate with
the general theory developed in the first two sections, one could
prove the coincidence of the virtual and actual thresholds of
canonical adjunction for varieties deviating from the
$K^2$-condition even stronger than those considered in this paper.
These limit resources of the techniques of the method of maximal
singularities will be used later.

\subsection {Historical remarks}

As soon as in [7] an effectively working technique of
investigating birational geometry of Fano fibrations $V/{\mathbb
P}^1$ satisfying the $K^2$-condition was developed, the
``boundary'' cases for which the deviation from the
$K^2$-condition was not too high, came up as a natural object of
further research. Already in [10] it was mentioned that for the
particular varieties which have already been studied the
$K^2$-condition was unnecessarily strong: the inequalities
ensuring birational rigidity have a considerable amount ``in
store''. Thus weakening this condition a little bit should not
change the final result. The papers of M.M.Grinenko [3-5] and
I.V.Sobolev [19,20] confirmed this idea. In [4] it was conjectured
that the $K$-condition (which is weaker than the $K^2$-condition)
is already sufficient for birational rigidity. This conjecture has
recently attracted new attention [1] in connection with the
attempts to study birational geometry of three-fold conic bundles
(over ${\mathbb P}^2$) which do not satisfy the Sarkisov condition
[17,18] (this condition is an exact analog of the $K^2$-condition
for conic bundles). However, in the papers [19,20] certain
varieties were successfully studied for which even the
$K$-condition was not true (although the deviation from this
condition was not too strong, either). And the technique used in
[19,20] was considerably weaker than that available today (see \S
2 of the present paper). This is an evidence that we do not
understand which mechanisms control birational rigidity (perhaps
one should speak of the ``degree of birational rigidity''). On the
other hand one can be optimistic concerning the prospects of
studying birational geometry of Fano fibrations by means of the
method of maximal singularities.

\subsection{Acknowledgements}

An essential part of this work was carried out at
Max-Planck-Institut f\" ur Mathematik in Bonn in the autumn 2003.
The author is very thankful to the Max-Planck-Institut f\" ur
Mathematik in Bonn for hospitality and the wonderful conditions of
work.


\section{The method of maximal singularities}

In this section we remind the main facts of the method of maximal
singularities: the Noether-Fano inequality, the concepts of a
maximal and a supermaximal singularities, the technique of
computing the self-intersection of a movable linear system
[8,13,15]. We prove Theorem 2 of the previous paper [16]. After
that, we modify the technique for the case when the
$K^2$-condition does not hold.

\subsection{Maximal singularities of linear systems}

Let $\Sigma\subset|-nK_V+lF|$, $l\in{\mathbb Z}_+$, be a movable
linear system on the variety $V$. Since $l\geq 0$, we get
$$
c(\Sigma)= n
$$
whereas $n=0$ if and only if the linear system $\Sigma$ comes from
the pencil $|F|$, that is, if it is pulled back from the base
${\mathbb P}^1$.

{\bf Remark 1.1.} If the fibration $V/{\mathbb P}^1$ satisfies the
$K^2$-condition, that is,
$$
K^2_V\not\in \mathop{\rm Int}A^2_+V,
$$
then for any movable system $\Sigma\subset|-nK_V+lF|$ we have
$l\geq 0$. Indeed, the self-intersection of the linear system
$\Sigma$
$$
(-nK_V+lF)^2=n^2K^2_V+2nlH_F
$$
is the class of an effective cycle of codimension two. By the
$K^2$-condition this implies that $l\in{\mathbb Z}_+$.

Assume that the inequality
$$
c_{virt}(\Sigma)< c(\Sigma)= n,
$$
holds, that is, there is a model $V^{\sharp}$ of the variety $V$
such that
$$
c(\Sigma^{\sharp},V^{\sharp})<c(\Sigma,V).
$$

{\bf Proposition 1.1.} {\it There is a prime divisor $E\subset
V^{\sharp}$, satisfying the Noether-Fano inequality}
\begin{equation}
\label{a1} \nu_E(\Sigma)>n\cdot a(V,E)
\end{equation}

The geometric discrete valuation $\nu_E$ of the field of rational
functions ${\mathbb C}(V)$, or any prime divisor $E^+\subset V^+$
on any model $V^+$ of the variety $V$, realizing this discrete
valuation, is called a {\it maximal singularity} of the linear
system $\Sigma$.

For a {\bf proof} of Proposition 1.1, see any of the papers
[6-8,13,15].

Note that the expression $a(V,E)$ in (\ref{a1}) denotes the
discrepancy of $E$ with respect to the original model $V$. Thus
the log-pair
$$
(V,\frac{1}{n}\Sigma)
$$
is not canonical: the exceptional divisor $E$ realizes its
singularity which is not canonical.

Since the linear system $\Sigma$ is movable and $\nu_E(\Sigma)>0$,
we conclude that the centre of the discrete valuation $\nu_E$ on
$V$ is a subvariety $B=\mathop{\rm centre}_V(\nu_E)\subset V$ of
codimension at least two. Let
\begin{equation}
\label{a2}
\begin{array}{rccc}
\varphi_{i,i-1}\colon & V_i & \to & V_{i-1}\\
  & \bigcup & & \bigcup \\
  & E_i & \to & B_{i-1}
\end{array}
\end{equation}
be the sequence of blow ups with irreducible centres
$B_{i-1}\subset V_{i-1}$, which is determined in a unique way by
the following conditions:

(1) $V_0=V$, $B_0=B$, $i=1,\dots,K$;

(2) $B_j=\mathop{\rm centre}_{V_j}(\nu_E)\subset V_j$,
$E_{j+1}=\varphi^{-1}_{j+1,j}(B_j)$;

(3) the valuation $\nu_{E_K}$ coincides with $\nu_E$.

\noindent In other words, the birational map
$$
V_K-\,-\,\to V^{\sharp}
$$
is biregular at the generic point of the divisor $E_K$ and
transforms $E_K$ into $E$. The symbol $\Sigma^j$ below means the
strict transform of the linear system $\Sigma$ on $V_j$. Set
$$
\nu_j= \mathop{\rm mult}\nolimits_{B_{j-1}}\Sigma^{j-1}, \quad
\delta_j= \mathop{\rm codim}B_{j-1}-1.
$$
On the set of exceptional divisors
$$
\{E_1,\dots,E_K\}
$$
we define in the usual way [6,8,13] an oriented graph structure:
an oriented edge (an arrow) goes from $E_i$ to $E_j$, if and only
if $i>j$ and
$$
B_{i-1}\subset E^{i-1}_j,
$$
which is denoted as $i\to j$. As usual, for $i>j$ set
$$
p_{ij}= \sharp\{\mbox{the paths from}\,\,  E_i\,\,\mbox{to}\,\,
E_j\}\geq 1,
$$
$p_{ii}=1$ by definition. Set $p_i=p_{Ki}$. The Noether-Fano
inequality takes the traditional form
$$
\sum^{K}_{i=1}p_i\nu_i>n\sum^{K}_{i=1}p_i\delta_i.
$$


\subsection{A stronger version of the Noether-Fano inequality}

Now assume that the general fiber $F=F_t$ of the fiber space
$V/{\mathbb P}^1$ admits no movable linear system with a maximal
singularity, that is, for any movable system
$\Sigma_F\subset|nH_F|=|-nK_F|$ and any geometric discrete
valuation $\nu_{E^*}$ the inequality
$$
\nu_{E^*}(\Sigma_F)\leq na(E^*)
$$
holds.

{\bf Proposition 1.2.} {\it The centre $B$ of the maximal
singularity $\nu_E$ on $V$ is contained in some fiber}
$\pi^{-1}(t)=F_t$, $t\in{\mathbb P}^1$.

{\bf Proof.} Assume the converse: $\pi(B)={\mathbb P}^1$. Let
$F\subset V$ be a fiber of general position. It is easy to see
that the restriction $\Sigma_F=\Sigma|_F$ of the linear system
$\Sigma$ onto $F$ is a movable linear system $\Sigma_F\subset
|nH_F|$ with a maximal singularity $\nu_{E^*}=\nu_E|_F$. The
simplest way to define this singularity is to restrict the
sequence of blow ups (\ref{a2}) onto the fiber $F$ and note that
the discrepancy remains the same:
$$
a(E|_F,F)=a(E,V)
$$
--- precisely for the reason that $B$ covers the base. The centre
of the valuation $\nu_E|_F$ is $B_F=B\cap F$. However, this
conclusion contradicts the assumption above. Therefore,
$\pi(B)\neq {\mathbb P}^1$. Q.E.D. for the proposition.

Let ${\cal M}=\{T_1,\dots,T_k\}$ be the set of all prime divisors
on $V^{\sharp}$, satisfying the Noether-Fano inequality (see
Proposition 1.1). As we have just proved, the centre
$B_E=\mathop{\rm centre}(\nu_E)$ of each maximal singularity
$E\in{\cal M}$ is contained in some fiber $F_t$. The set ${\cal
M}$ is finite (since the model $V^{\sharp}$ is fixed), so that
there is at most finite set of points $t\in{\mathbb P}^1$, the
fibers $F_t$ over which contain the centres of maximal
singularities. Set ${\cal M}_t=\{E\in{\cal M}|B_E\subset F_t\}$,
$$
e(E)=\nu_E(\Sigma)-na(E,V)>0
$$
for $E\in{\cal M}$. Recall that $\Sigma\subset|-nK_V+lF|$,
$l\in{\mathbb Z}_+$.

{\bf Proposition 1.3.} {\it The following inequality holds:}
\begin{equation}
\label{a3} \sum_{t\in{\mathbb P}^1}\mathop{\rm max}_{ \{E\in{\cal
M}_t\} } \frac{e(E)}{\nu_E(F_t)}>l
\end{equation}

{\bf Proof.} Let $D^{\sharp}\in\Sigma^{\sharp}$ be a general
divisor, that is, the strict transform on $V^{\sharp}$ of a
divisor $D\in\Sigma$ of general position. By assumption, the
linear system
$$
|D^{\sharp}+nK_{V^{\sharp}}|
$$
is empty. Therefore the linear system
$$
|lF-\sum_{E\in{\cal M}}e(E)E|
$$
is empty, too. On the other hand, by construction for $E\in{\cal
M}_t$ the divisor
$$
F_t-\nu_E(F_t)E
$$
is effective, so that the divisor
$$
\sum_{t\in{\mathbb P}_1}[\left( \mathop{\rm max}_{\{E\in{\cal
M}_t\}}\frac{e(E)}{\nu_E(F_T)}\right)F_t-\sum_{E\in{\cal
M}_t}e(E)E]
$$
is also effective. This immediately implies the inequality
(\ref{a3}). Q.E.D. for the proposition.

\subsection{The self-intersection of the linear system $\Sigma$}

Let $D_i\in\Sigma$, $i=1,2$, be general divisors, so that the
closed set $D_1\cap D_2$ is of codimension two. Let
$$
Z=(D_1\circ D_2)=Z^v+Z^h
$$
be the decomposition of the algebraic cycle of the
scheme-theoretic intersection of these divisors into the vertical
$(Z^v)$ and horizontal $(Z^h)$ parts. For the cycle $Z^v$ we have
a further decomposition
$$
Z^v=\sum_{t\in{\mathbb P}^1}Z^v_t,\quad \mathop{\rm
Supp}Z^v_t\subset F_t.
$$
Let $E\in{\cal M}_t\subset{\cal M}$ be a maximal singularity. Fix
$t$ and $E$ and apply the {\it technique of counting
multiplicities}, developed in [8,10,13], to the cycle $Z^v_t+Z^h$.
We assume that $t$ and $E$ are fixed throughout this subsection,
so that we write below $F$, $Z^v$, $e$, $B$ instead of $F_t$,
$Z^v_t$, $e(E)$, $B_E=\mathop{\rm centre}(E)$, respectively.

{\bf Lemma 1.1.} {\it The following estimate holds:}
$$
\mathop{\rm codim}\nolimits_{F} B\geq 2
$$

{\bf Proof.} Assume the converse: $B\subset F$ is a prime divisor.
Let $D\in \Sigma$ be a general divisor, $D_F$ its restriction onto
$F$. By the Noether-Fano inequality $\mathop{\rm mult}_BD>n$, so
that
$$
D_F=\alpha B+D^{\sharp},
$$
where $\alpha>n$ and $D^{\sharp}$ is an effective divisor on $F$.
However,
$$
D_F\sim nH_F,
$$
which gives an immediate contradiction. Q.E.D. for the lemma.

Now consider the sequence of blow ups (\ref{a2}) associated with
the discrete valuation $E$. We use the notations of Sec. 1.1. The
strict transforms $(Z^h)^j$, $(Z^v)^j$ of the cycles $Z^h$, $Z^v$
on $V_j$ are well defined. By the symbol $F^j$ we denote the
strict transform of the fiber $F$ on $V_j$. Set
$$
N=\mathop{\rm max}\{i\,|\,B_{i-1}\subset F^{i-1}\}.
$$
It is easy to see that $\varphi_{i,i-1}(B_i)=B_{i-1}$ for any
$i=1,\dots,K-1$, so that the codimensions $\mathop{\rm codim}B_i$
do not increase. Set
$$
L=\mathop{\rm max}\{i\,|\,\mathop{\rm codim}B_{i-1}\geq 3\}\leq K.
$$
We use also the following notations: for $i\in\{1,\dots,L\}$
$$
m^h_i=\mathop{\rm mult}\nolimits_{B_{i-1}}(Z^h)^{i-1}, \quad
m^v_i=\mathop{\rm mult}\nolimits_{B_{i-1}} (Z^v)^{i-1},
$$
where $m^{h(v)}_i\leq m^{h(v)}_{i-1}$ for $i=2,\dots,L$.

Set also
$$
\mu_i=\mathop{\rm mult}\nolimits_{B_{i-1}}F^{i-1}.
$$
Obviously, $\mu_i=0$ for $i\geq N+1$. The more so, $m^v_i=0$ for
$i\geq N+1$ (if $N<L$). Using the symbol $p_i$, as above, for the
number of paths in the graph $\Gamma$ of the resolution (\ref{a2})
from the vertex $E=E_K$ to $E_i$, we get

{\bf Proposition 1.4.} {\it The following inequality holds:}
\begin{equation}
\label{a4}
\sum^L_{i=1}p_im^h_i+\sum^{\min\{N,L\}}_{i=1}p_im^v_i\geq
\sum^K_{i=1}p_i\nu^2_i\geq
\frac{\displaystyle(n\sum^K_{i=1}p_i\delta_i+e)^2}{\displaystyle\sum^K_{i=1}p_i}
\end{equation}

{\bf Proof} is obtained by applying the technique of counting
multiplicities [8,10,13] combined with the condition
$$
\sum^K_{i=1}p_i\nu_i=n\sum^K_{i=1}p_i\delta_i+e,
$$
where $e>0$. Here we do not repeat these standard arguments, just
referring to [10].

\subsection{$K^2$-condition and birational rigidity}

As a first example of using the technique described above let us
prove Theorem 2 of the previous paper [16]. The arguments below
follow the lines of the proof of particular cases of this theorem,
given in [7,10] for certain special families of Fano fiber spaces.
Here we consider the general case.

Assume that the fiber space $V/{\mathbb P}^1$ satisfies the
$K^2$-condition: $K^2_V\not\in\mathop{\rm Int}A^2_+V$. Then for
some $\alpha\in{\mathbb Z}_+$ the following relation holds:
$$
Z^h\sim n^2K^2_V+\alpha H_F,
$$
so that
\begin{equation}
\label{a5} \mathop{\rm deg}Z^v=\sum_{t\in{\mathbb P}^1}\mathop{\rm
deg}Z^v_t\leq(2n\mathop{\rm deg}V)l.
\end{equation}

{\bf Proposition 1.5.} {\it For some point $t\in{\mathbb P}^1$
there is a maximal singularity $E\in{\cal M}_t$, satisfying the
estimate}
\begin{equation}
\label{a6} e(E)>\frac{\nu_E(F_t)}{2n\mathop{\rm deg}V}\mathop{\rm
deg}Z^v_t.
\end{equation}

Following [7,10], we call the singularity $E$, satisfying the
inequality (\ref{a6}), a {\it supermaximal singularity}.

{\bf Proof of Proposition 1.5:}  it is sufficient to compare the
inequalities (\ref{a3}) and (\ref{a5}). Q.E.D.

Now in the notations of subsection 1.3 set
$$
\Sigma_{l}=\sum^L_{i=1}p_i,\quad
\Sigma_{u}=\sum^K_{i=L+1}p_i,\quad
\Sigma_{f}=\sum^{\min\{N,L\}}_{i=2}p_i.
$$
Note that
$$
\nu_E(F)=\sum^N_{i=1}p_i\mu_i\leq p_1\mu_1+\mu_2\Sigma_{f}
$$
by definition of the multiplicities $\mu_i$. Obviously,
$$
m^h_i\leq m_h=m^h_1=\mathop{\rm mult}\nolimits_{B}Z^h.
$$
Set also
$$
d_h=\mathop{\rm deg}Z^h,\quad d_v=\mathop{\rm deg}Z^v_t
$$
and introduce the coefficients
$$
k_h=\frac{m_h}{d_h}\mathop{\rm deg} V,\quad k_v=\frac{\mathop{\rm
deg}V}{\nu_E(F)d_v}\sum^{\min\{N,L\}}_{i=1}p_im^v_i.
$$
Now Proposition 1.5 implies

{\bf Corollary 1.1.} {\it The following estimate holds:}
$$
(4-k_h)\Sigma_l(\Sigma_l+\Sigma_u)n^2+\Sigma^2_u
n+e^2+2(2-k_v)\Sigma_lne+2(1-k_v)\Sigma_une<0.
$$

{\bf Proof.} In the inequality (\ref{a4}) replace $m^h_i$ by
$m_h=k_hd_h/\mathop{\rm deg}V=k_hn^2$, the numbers $\delta_i$ for
$i\leq L$ by 2, which could only make the inequality sharper. Now
taking into account the definition of the coefficient $k_v$ after
easy computations we get
\begin{equation}
\label{a7}
\begin{array}{c}
\displaystyle (4-k_h)\Sigma_l(\Sigma_l+\Sigma_u)n^2+(n\Sigma_u+e)^2+4\Sigma_len-\\
   \\
\displaystyle - k_v\frac{d_v\nu_E(F)}{\mathop{\rm
deg}V}(\Sigma_l+\Sigma_u)\leq 0.
\end{array}
\end{equation}

Taking into consideration the definition of a supermaximal
singularity (Proposition 1.5), replace $d_v\nu_E(F)$ by
$2ne\mathop{\rm deg}V$. This makes our inequality a strict one and
we obtain exactly what we claimed. Q.E.D. for the corollary.

Now assume that the fiber space $V/{\mathbb P}^1$ satisfies

\begin{itemize}
\item the condition (v) of the paper [16], that is, for any
irreducible vertical subvariety $Y$ of codimension 2, $Y\subset
\pi^{-1}(t)=F_t$, and any smooth point $o\in F_t$ the estimate
$$
\frac{\mathop{\rm mult}\nolimits_o}{\mathop{\rm deg}} Y \leq
\frac{2}{\mathop{\rm deg}\nolimits V}
$$
holds;

\item the condition (vs) of the paper [16], that is, for any
vertical subvariety $Y\subset F_t$ of codimension 2 (with respect
to $V$, that is, for a prime divisor on $F_t$), a singular point
$o\in F_t$ and an infinitely near point $x\in \widetilde F_t$,
where $\varphi\colon \widetilde F_t\to F_t$ is a blow up of the
point $o$, $\varphi(x)=o$, $\widetilde Y\subset \widetilde F_t$
the strict transform of the subvariety $Y$ on $\widetilde F_t$,
the following estimates hold:
$$
\frac{\mathop{\rm mult}\nolimits_o}{\mathop{\rm deg}} Y \leq
\frac{4}{\mathop{\rm deg}\nolimits V}, \quad \frac{\mathop{\rm
mult}\nolimits_x \widetilde Y }{\mathop{\rm deg}\nolimits Y} \leq
\frac{2}{\mathop{\rm deg}\nolimits V};
$$

\item the condition (h) of the paper [16], that is, for any
horizontal subvariety $Y$ of codimension 2 and a point $o\in Y$
the estimate
$$
\frac{\mathop{\rm mult}\nolimits_o}{\mathop{\rm deg}} Y \leq
\frac{4}{\mathop{\rm deg}\nolimits V}.
$$
holds.

\end{itemize}

In these assumptions we get

{\bf Lemma 1.2.} {\it The following estimates hold:} $k_h\leq 4$,
$k_v\leq 2$.

{\bf Proof.} The first inequality follows directly from the
definition of the number $k_h$ and the condition (h), which is by
assumption satisfied for the fiber space $V/{\mathbb P}^1$. Let us
prove the second inequality. We get
$$
\left(\frac{d_v}{\mathop{\rm deg}V}\right)
k_v=\frac{\displaystyle\sum^{\mathop{\rm min}
\{N,L\}}_{i=1}p_im^v_i}{\displaystyle\sum^N_{i=1}p_i\mu_i}\leq
\frac{p_1m^v_1+\Sigma_{f}m^v_2} {p_1\mu_1+\Sigma_{f}}
$$

If $\mathop{\rm dim}B\geq 1$ or $B=o\in F$ is a non-singular point
of the fiber, then $\mu_1=\dots=\mu_N=1$ and by the condition (v)
$$
m^v_2\leq m^v_1\leq 2\frac{d_v}{\mathop{\rm deg}V},
$$
which immediately implies the inequality $k_v\leq 2$. If $B=o\in
F$ is a singular point of the fiber, then $\mu_1\geq 2$ and by the
condition (vs)
$$
m^v_1\leq \frac{4d_v}{\mathop{\rm deg}V},\quad m^v_2\leq
\frac{2d_v}{\mathop{\rm deg}V},
$$
whence we get again that $k_v\leq 2$. Q.E.D. for the lemma.

Recall the claim of Theorem 2 of the paper [16]:

{\it Assume that a smooth standard Fano fiber space $V/{\mathbb
P}^1$ satisfies the $K^2$-condition and the conditions {\rm (v)},
{\rm (vs)} and {\rm (h)}. Then $V/{\mathbb P}^1$ is birationally
superrigid.}

Let us complete the proof of this theorem. Assume the converse.
Then Corollary 1.1 and Lemma 1.2 give us the inequality
$$
\Sigma^2_un^2-2\Sigma_une+e^2<0.
$$
This is, however, impossible. We get a contradiction. Q.E.D. for
Theorem 2 of the paper [16].

\subsection{The generalized $K^2$-condition}

The methods developed in [7,10] and reproduced above in the
general form work well for Fano fiber spaces $V/{\mathbb P}^1$
which do not satisfy the $K^2$-condition. If the deviation from
the $K^2$-condition is not too great then the technique of Sec.
1.4 is still effective and requires only a slight modification.
First of all, one should be able to measure the deviation from the
$K^2$-condition.

{\bf Definition 1.1.} A standard Fano fiber space $V/{\mathbb
P}^1$ satisfies the {\it generalized $K^2$-condition of depth}
$\varepsilon\geq 0$, if
$$
K^2_V-\varepsilon H_F\not\in \mathop{\rm Int}A^2_+V.
$$

As mentioned in [10], for effective cycles of codimension 2 on the
natural Fano varieties sharper estimates on the multiplicities are
satisfied than one needs to prove birational rigidity. However,
for Fano fiber spaces that do not satisfy the $K^2$-condition
these estimates turn out to be useful.

{\bf Definition 1.2.} A standard Fano fibration $V/{\mathbb P}^1$
satisfies the {\it generalized condition} (h) {\it of depth}
$\delta\geq 0$, if for any horizontal subvariety $Y\subset V$ of
codimension two and an arbitrary point $o\in Y$ the following
inequality holds:
$$ \frac{\mathop{\rm
mult}\nolimits_o}{\mathop{\rm deg}}Y\leq
\frac{4-\delta}{\mathop{\rm deg}V}.
$$

Starting from this moment and throughout this section we assume
that the fiber space $V/{\mathbb P}^1$ satisfies the generalized
$K^2$-condition of depth $\varepsilon\geq 0$ and the generalized
condition (h) of depth $\delta\geq 0$. Besides, we assume that the
conditions (v) and (vs) are satisfied (in their normal form). Fix
a movable linear system $\Sigma\subset |-nK_V+lF|$ with
$l\in{\mathbb Z}_+$ and assume that $c_{\rm virt}(\Sigma)<n$. For
the horizontal part of the self-intersection of the linear system
$\Sigma$ we get
$$
Z^h\sim n^2K^2_V+\alpha H_F,
$$
where the coefficient $\alpha\in{\mathbb Z}$ satisfies the
inequality
$$
\alpha\geq -\varepsilon n^2.
$$
Therefore, for the vertical component we get
$$
Z^v\sim (2nl-\alpha)H_F,
$$
whence
\begin{equation}
\label{a8} \mathop{\rm deg}Z^v=\sum_{t\in{\mathbb P}^1}\mathop{\rm
deg}Z^v_t\leq(2nl+\varepsilon n^2)\mathop{\rm deg}V.
\end{equation}

{\bf Proposition 1.6.} {\it For some point $t\in{\mathbb P}^1$
there exists a maximal singularity $E\in{\cal M}_t\neq \emptyset$,
satisfying the estimate}
\begin{equation}
\label{a9} e(E)>\frac{\nu_E(F_t)}{2}\left(\frac{\mathop{\rm
deg}Z^v_t}{n\mathop{\rm deg}V}-\varepsilon n\right)
\end{equation}

{\bf Remark 1.2.} For $\varepsilon=0$ we get Proposition 1.5.

{\bf Proof of Proposition 1.6.} Compare the inequalities
(\ref{a3}) and (\ref{a8}). Replacing the number $l$ in the
right-hand side of the inequality (\ref{a8}) by the left-hand side
of the inequality (\ref{a3}), we get
$$
\sum_{t\in{\mathbb P}^1}\left[\mathop{\rm deg} Z^v_t-2n\mathop{\rm
deg}V\mathop{\rm max}_{\{E\in{\cal M}_t\}}\frac{e(E)}{\nu_E(F_t)}
\right]< \varepsilon n^2\mathop{\rm deg}V,
$$
which immediately implies our claim. Q.E.D. for Proposition 1.6.

{\bf Remark 1.3.} If there are a few maximal singularities the
centres of which lie in the fibers over distinct points
$t_1,\dots,t_k$ then the claim of Proposition 1.6 can be improved:
there is a maximal singularity $E\in{\cal M}_t$,
$t\in\{t_1,\dots,t_k\}$ satisfying the estimate
$$
e(E)>\frac{\nu_E(F_t)}{2}\left(\frac{\mathop{\rm
deg}Z^v_t}{n\mathop{\rm deg}V}-\frac{\varepsilon n}{k}\right).
$$
However, this improvement is hardly useful, since to prove
birational rigidity, that is, to realize the full scheme of the
method of maximal singularities, the worst case should be
considered.

Now we argue as in the proof of Corollary 1.1 with the only
difference: the expression $d_v\nu_E(F)$ in the inequality
(\ref{a7}) is replaced by the expression
$$
(2ne+\varepsilon n^2\nu_E(F))\mathop{\rm deg}V,
$$
which makes the inequality sharp. Taking into account that the
conditions (h), (v) and (vs) still hold, so that the claim of
Lemma 1.2 is true, we get the inequality
$$
((4-k_h)\Sigma_l-k_v\varepsilon\nu_E(F))(\Sigma_l
+\Sigma_u)n^2+(n\Sigma_u-e)^2<0.
$$
By the definition of the number $k_v$ we get
$$
k_v\nu_E(F)=\mathop{\rm
deg}V\sum^{\min\{N,L\}}_{i=1}p_i\frac{m^v_i}{d_v}\leq2\Sigma_l.
$$
Therefore, the following inequality holds:
\begin{equation}
\label{a10}
(4-k_h-2\varepsilon)\Sigma_l(\Sigma_l+\Sigma_u)n^2+(n\Sigma_u-e)^2<0.
\end{equation}
Now recall that by the generalized condition (h) of depth
$\delta\geq 0$ the coefficient $k_h$ satisfies the estimate
$$
k_h\leq 4-\delta.
$$
This immediately implies

{\bf Proposition 1.7.} {\it If $\delta\geq 2\varepsilon$, then the
equality
$$
c_{\rm virt}(\Sigma)=c(\Sigma)=n
$$
holds. In particular, if for any movable linear system
$\Sigma\subset |-nK_V+lF|$ we have $l\in{\mathbb Z}_+$, then the
fiber space $V/{\mathbb P}^1$ is birationally superrigid.}


\section{The technique of counting multiplicities}

The aim of this section is to sharpen the inequality (\ref{a10})
and the claim of Proposition 1.7. For this purpose we take into
consideration the dimensions of the blown up subvarieties
$B_{i-1}$: the smaller is the dimension of the centre of a blow
up, the higher is the discrepancy of the exceptional divisor and,
therefore, the better is the estimate for the multiplicity
$\nu_E(\Sigma)$. On the other hand, the only working method of
getting an upper bound for the singularities of the horizontal
cycle $Z^h$ is to restrict $Z^h$ onto the fiber $F$ and estimate
the singularities of the effective cycle $(Z^h\circ F)$ of
codimension two on $F$. This method was used above; also in the
papers [7,10,16] the singularities of the cycle $Z^h$ were
estimated in this way. However, up to this day we never took into
account the input of the subvarieties $B_{i-1}$ of codimension
three (with respect to $V$), lying in the strict transform
$F^{i-1}$ of the fiber $F$ (provided they exist). That is what we
do below. As a result, we obtain estimates which make it possible
to exclude a maximal singularity even in cases when the deviation
from the $K^2$-condition is essential: it is sufficient that the
generalized $K^2$-condition of depth 2 holds.


\subsection{The notations and the principal claim}

We go on studying the movable linear system $\Sigma\subset
|-nK_V+lF|$, $l\in{\mathbb Z}_+$, satisfying the inequality
$$
c_{\rm virt}(\Sigma)<c(\Sigma)=n.
$$
All notations of \S 1 are valid. However, if in subsections 1.4
and 1.5 we treated all blow ups of subvarieties of codimension
three and higher in the same way, now we argue in a more refined
way. Set
$$
J_s=\{i\,|\,1\leq i\leq K,\quad \mathop{\rm codim}B_{i-1}\geq 4\},
$$
$$
J_m=\{i\,|\,1\leq i\leq K,\quad \mathop{\rm codim}B_{i-1}=3\},
$$
$$
J_u=\{i\,|\,L+1\leq i\leq K\},\quad J_l=J_s\cup J_m.
$$
In its turn, let us break the set $J_m$ into two disjoint subsets,
$J_m=J^+_m \coprod J^-_m$, where
$$
J^+_m=\{i\in J_m\,|\,B_{i-1}\subset F^{i-1}\},
$$
$J^-_m=J_m\setminus J^+_m=\{i\in J_m\,|\,B_{i-1}\not\subset
F^{i-1}\}$. It might well turn out that the set $J^+_m$ or $J^-_m$
(or the whole set $J_m$) is empty. Set furthermore
$$
\Sigma_s=\sum\limits_{i\in J_s}p_i,\quad
\Sigma^{\pm}_m=\sum\limits_{i\in J^{\pm}_m}p_i,\quad
\Sigma_m=\Sigma^+_m+\Sigma^-_m,
$$
whereas the symbol $\Sigma_u$ retains its previous meaning. In the
notations of Sec. 1.4 we get $\Sigma_l=\Sigma_s+\Sigma_m$. Now the
inequality (\ref{a4}) can be rewritten as
\begin{equation}\label{b1}
\sum_{i\in J_l}p_im^h_i+\sum_{i\in J_S\cup J^+_m}p_im^v_i\geq
\frac{((3\Sigma_s+2\Sigma_m+\Sigma_u)n+e)^2}{\Sigma_s+\Sigma_m+\Sigma_u}
\end{equation}
Recall that $\mu_i=\mathop{\rm mult}\nolimits_{B_{i-1}}F^{i-1}$,
where $\mu_i=1$ for $i\geq 2$ and for $\mu_1$ there are two
possible cases: $\mu_1=1$ or $\mu_1=2$.

{\bf Proposition 2.1.}  {\it The following estimate holds:
\begin{equation}\label{b2}
\sum_{i\in J_s\cup J^+_m}\mu_ip_im^h_i\leq p_1\mathop{\rm
mult}\nolimits_B(Z^h\circ F)+(\Sigma_s-p_1)\mathop{\rm
mult}\nolimits_{B_1}(Z^h\circ F)^1.
\end{equation}
In particular,}
\begin{equation}\label{b3}
\sum_{i\in J_s\cup J^+_m}\mu_ip_im^h_i\leq \Sigma_s\mathop{\rm
mult}\nolimits_B(Z^h\circ F).
\end{equation}


\subsection{Proof of Proposition 2.1: counting multiplicities}

To begin with, let us consider the following general situation.
Let $Y\subset V$ be an irreducible horizontal subvariety of
codimension two, $Y^i\subset V_i$ its strict transform,
\begin{equation}\label{b4}
m_Y(i)=\mathop{\rm mult}\nolimits_{B_{i-1}}Y^{i-1}
\end{equation}
the corresponding multiplicity. Set $Y_F=(Y\circ F)$. This is an
effective class of codimension two in the fiber $F$. Let
$Y^i_F\subset V_i$ be its strict transform and
\begin{equation}\label{b5}
m_{Y,F}(i)=\mathop{\rm mult}\nolimits_{B_{i-1}}Y^{i-1}.
\end{equation}
Since the support of the cycle $Y_F$ is contained in the fiber
$F$, the numbers $m_{Y,F}(i)$ vanish for $i\in J^-_m$.

{\bf Lemma 2.1.} {\it The following estimate holds:}
\begin{equation}\label{b6}
\sum_{i\in J_s\cup  J^+_m}p_im_Y(i)\mu_i\leq \sum_{i\in
J_s}p_im_{Y,F}(i).
\end{equation}

Before we start to prove it, recall some facts which follow
immediately from elementary intersection theory [21]. There are no
convenient references because here we consider the case when a
divisor intersects a subvariety, whereas in [8,13] the case of two
divisors was treated. Let $X$ be an arbitrary smooth variety,
$B\subset X$, $B\not\subset\mathop{\rm Sing} X$, an irreducible
subvariety of codimension $\geq 2$, $\sigma_B\colon X(B)\to X$ its
blow up, $E(B)=\sigma^{-1}_B(B)$ the exceptional divisor. Let
$$
Z=\sum m_iZ_i, \quad Z_i\subset E(B),
$$
be a cycle of dimension $k$, $k\geq\mathop{\rm dim} B$. We define
the {\it degree} of the cycle $Z$, setting
$$
\mathop{\rm deg} Z= \sum_im_i\mathop{\rm deg}\left(
Z_i\bigcap\sigma^{-1}_B(b) \right),
$$
where $b\in B$ is a point of general position,
$\sigma^{-1}_B(b)\cong{\mathbb P}^{\mathop{\rm codim} B-1}$ and
the degree in the right-hand side is the usual degree in the
projective space.

Note that $\mathop{\rm deg} Z_i=0$ if and only if $\sigma_B(Z_i)$
is a proper closed subset of the subvariety $B$.

Now let $D$ be a prime Weil divisor on $X$, $Y\subset X$ an
irreducible subvariety of dimension $l\leq \mathop{\rm dim}X-1$.
Assume that $Y\not\subset D$ and that $\mathop{\rm dim}B\leq l-1$.
The strict transforms of the divisor $D$ and the subvariety $Y$ on
$X(B)$ are denoted by the symbols $D^B$ and $Y^B$, respectively.

{\bf Lemma 2.2.} (i) {\it Assume that $\mathop{\rm dim}B\leq l-2$.
Then
$$
D^B\circ Y^B=(D\circ Y)^B+Z,
$$
where $\circ$ means the operation of taking the algebraic cycle of
the scheme-theoretic intersection, $\mathop{\rm Supp}Z\subset
E(B)$ and
$$
\mathop{\rm mult}\nolimits_B(D\circ Y)=\mathop{\rm
mult}\nolimits_BD\cdot\mathop{\rm mult}\nolimits_BY+\mathop{\rm
deg}Z.
$$

{\rm (ii)} Assume that $\mathop{\rm dim}B=l-1$. Then
$$
D^B\circ Y^B=Z+Z_1,
$$
where $\mathop{\rm Supp}Z\subset E(B)$, $\mathop{\rm
Supp}\sigma_B(Z_1)$ does not contain $B$ and
$$
D\circ Y= \left[ (\mathop{\rm mult}\nolimits_BD)(\mathop{\rm
mult}\nolimits_BY)+ \mathop{\rm deg} Z \right] B+ (\sigma_B)_*Z_1.
$$
}\vspace{0.1cm}

{\bf Proof} is easily obtained by means of the standard
intersection theory  [21].


\subsection{Proof of Lemma 2.1}

Let us construct a sequence of effective cycles of codimension
three on the varieties $V_i$, setting
$$
\begin{array}{ccl}
\displaystyle
Y\circ F & = & Z_0\,\, (=Y_F),\\
\displaystyle
Y^1\circ F^1 & = & Z^1_0+Z_1,\\
\displaystyle
\vdots\\
\displaystyle Y^i\circ F^i & = &
(Y^{i-1}\circ F^{i-1})^i+Z_i,\\
\displaystyle \vdots
\end{array}
$$
$i\in J_s$, where $\mathop{\rm Supp} Z_i\subset E_i$. Thus for any
$i\in J_s$ we get:
$$
Y^i\circ F^i= Y^i_F+Z^i_1+\dots+Z^i_{i-1}+Z_i.
$$
For any $j>i$, $j\in J_s$ set
$$
m_{i,j}=\mathop{\rm mult}\nolimits_{B_{j-1}}(Z^{j-1}_i)
$$
(the multiplicity of an irreducible subvariety along a smaller
subvariety is understood in the usual sense; for an arbitrary
cycle we extend the multiplicity by linearity).

Now set $d_i=\mathop{\rm deg} Z_i$. We get the following system of
equalities:
$$
\begin{array}{l}
\displaystyle
m_Y(1)\mu_1+ d_1= m_{Y,F}(1),\\
\displaystyle
m_Y(2)\mu_2+d_2=m_{Y,F}(2)+m_{1,2},\\
\displaystyle
\vdots\\
\displaystyle
m_Y(i)\mu_i+d_i=m_{Y,F}(i)+m_{1,i}+\dots+m_{i-1,i}\\
\displaystyle
\vdots\\
\end{array}
$$
for all $i\in J_s$. Setting $S=\max\{i\in J_s\}$, look at the last
equality in this sequence:
$$
m_Y(S)\mu_S+d_S=m_{Y,F}(S)+m_{1,S}+\dots+m_{S-1,S}.
$$
If $J^+_m\neq \emptyset$, then, by the part (ii) of Lemma 2.2, we
obtain
$$
d_S\geq\sum_{i\in J^+_m}m_Y(i)\mu_i\mathop{\rm
deg}(\varphi_{i-1,S})_*B_{i-1}\geq \sum_{i\in J^+_m}m_Y(i)\mu_i.
$$
Recall the following useful

{\bf Definition 2.1.} (see [8,13]). A function $a\colon J_s\to
{\mathbb R}_+$ is said to be {\it compatible with the graph
structure}, if
$$
 a(i)\geq\sum_{\scriptstyle \begin{array}{c}\scriptstyle
 j\to i,\\ \scriptstyle
j\in J_s
\end{array}}a(j)
$$
for any $i\in J_s$.

We will actually use only one function, compatible with the graph
structure, namely $a(i)=p_i$.

{\bf Proposition 2.2.} {\it Let $a(\cdot)$ be a function,
compatible with the graph structure. Then the following inequality
holds:}
\begin{equation} \label{b8}
\sum_{i\in J_s}a(i)m_{Y,F}(i)\geq \sum_{i\in
J_s}a(i)m_Y(i)\mu_i+a(S)\sum_{i\in J^+_m}m_Y(i)\mu_i.
\end{equation}

{\bf Proof} is obtained in exactly the same way as in the case of
two divisors ([8,13]) multiply the $i$-th equality by $a(i)$ put
them all together. In the right hand side for any $i\geq 1$ we
obtain the expression
$$
\sum_{j\geq i+1}a(j)m_{i,j}.
$$
In the left hand side for any $i\geq 1$ we obtain the component
$a(i)d_i$.

{\bf Lemma 2.3.} {\it If $m_{i,j}> 0$, then $j\to i$}.

{\bf Proof} [8,13]: if $m_{i,j}> 0$, then
$B_{j-1}\subset\mathop{\rm Supp}Z^{j-1}_i$, but $\mathop{\rm
Supp}Z_i\subset E_i$ so that $B_{j-1}\subset E^{j-1}_i$. Q.E.D.
for the lemma.

The next standard step is to compare the multiplicities $m_{i,j}$
with the degrees.

{\bf Lemma 2.4.} {\it For any} $i<j\in J_s$ {\it we get}
$$
m_{i,j}\leq d_i.
$$

{\bf Proof.} If $m_{i,j}=0$, then there is nothing to prove.
Otherwise $j\to i$ and we need to prove that
$$
\mathop{\rm mult}\nolimits_{B_{j-1}}Z^{j-1}_i\leq \mathop{\rm
deg}Z_i.
$$
Taking into account that the maps $\varphi_{a,b}\colon B_a\to B_b$
are surjective, it suffices to prove the inequality
\begin{equation}\label{b7}
\mathop{\rm
mult}\nolimits_{[B_{j-1}\cap\varphi^{-1}_{i,i-1}(t)^{j-1}]}
[Z_i\cap\varphi^{-1}_{i,i-1}(t)]^{j-1}\leq \mathop{\rm deg}
[Z_i\cap\varphi^{-1}_{i,i-1}(t)],
\end{equation}
where $t\in B_{i-1}$ is a point of general position. Taking into
account that $\varphi^{-1}_{i,i-1}(t)$ is the projective space
${\mathbb P}^{\mathop{\rm codim}B_{i-1}-1}$, we see that the right
hand side in (\ref{b7}) means the usual degree of a hypersurface
in the projective space whereas the set $[Z_i\cap
\varphi^{-1}_{i,i-1}(t)]^{j-1}$ is obtained from this hypersurface
by means of a finite sequence of blow ups
$\varphi_{s,s-1},s=i+1,\dots,j-1$, restricted onto
$\varphi^{-1}_{i,i-1}(t)$. Taking into consideration that the
multiplicities do not increase when blowing ups are performed, we
reduce the claim to the obvious case of a hypersurface in the
projective space. Q.E.D. for the lemma.

As a result, we obtain the following estimate
$$
\sum_{j\geq i+1}a(j)m_{i,j}= \sum_{\scriptstyle
\begin{array}{c}
\scriptstyle
j\geq i+1\\
\scriptstyle m_{i,j}\neq 0
\end{array}
}a(j)m_{i,j}\leq d_i\sum_{j\to i}a(j)\leq a(i)d_i.
$$
By what was said above, we may remove from the right hand side all
the components $m_{i,*},i\geq 1$ and from the left hand side all
the components $d_i, i\geq 1$, replacing the equality sign $=$ by
the inequality sign $\leq$. Q.E.D. \vspace{0.3cm}

Setting in the inequality (\ref{b8}) $a(i)=p_i$ and recalling that
for $j\geq S$ we have $p_j\leq p_S$, we complete the proof of
Lemma 2.1.

Now let us complete the proof of Proposition 2.1.

It is obvious that the inequality (\ref{b6}) remains valid if $Y$
is an effective horizontal cycle of codimension 2 on $V$, that is,
each component of the cycle $Y$ is a horizontal subvariety.
Besides, the formulae (\ref{b4}), (\ref{b5}) extend by linearity
to the set of all effective horizontal cycles, and the left hand
and the right hand sides of the inequality (\ref{b6}) are linear
in $m_Y(\cdot),m_{Y,F}(\cdot)$, respectively.

Now set $Y=Z^h$ and take into account that
$$
m_{Y,F}(i)\leq \mathop{\rm mult}\nolimits_{B_1}(Z^h\circ F)^1
$$
for $i\geq 2$. This implies the inequality (\ref{b2}). The second
inequality of Proposition 2.1 follows from (\ref{b2}). Q.E.D. for
the proposition.

{\bf Remark 2.1.} The inequality (\ref{b3}) is more compact than
(\ref{b2}), however in some cases it is possible to get a stronger
estimate for $\mathop{\rm mult}\nolimits_{B_1}(Z^h\circ F)^1$ than
for $\mathop{\rm mult}\nolimits_B(Z^h\circ F)$.


\subsection{Estimating the multiplicities of a linear system: \\ the
non-singular case}

Now let us turn to the main problem, that is, estimating the
singularities of the movable linear system $\Sigma$. Assume that
the fiber space $V/{\mathbb P}^1$ satisfies the conditions (v) and
(vs), and at any point $o\in V$ satisfies at least one of the
conditions (f) or (fs) formulated below:

(f) for any irreducible subvariety $Y\ni o\in F$ of codimension 2
(with respect to the fiber $F$) the following inequality holds:
\begin{equation}
\label{b9} \frac{\mathop{\rm mult}\nolimits_o}{\mathop{\rm
deg}}Y\leq\frac{4}{\mathop{\rm deg}V},
\end{equation}

(fs) $o\in F$ is a double point of the fiber and for any
irreducible subvariety $Y\subset F$  of codimension 2 the
following estimate holds:
\begin{equation}
\label{b14} \frac{ \mathop{\rm mult}\nolimits_o }{ \mathop{\rm
deg}}Y\leq \frac{6}{\mathop{\rm deg}V},
\end{equation}
and, besides, for any infinitely near point of the first order
$x\in E_F$, where $\varphi_o\colon \widetilde F\to F$ is the blow
up of the point $o\in F$, $E_F \subset \widetilde F$ the
exceptional divisor, the inequality
$$
\frac{\mathop{\rm mult}\nolimits_x\widetilde Y}{\mathop{\rm
deg}Y}\leq \frac{3}{\mathop{\rm deg}V}
$$
holds, where $\widetilde Y$ is the strict transform of the
subvariety $Y$ on $\widetilde F$.

{\bf Remark 2.2.} Let us draw the reader's attention to the fact
that in the condition (f) we do not specify whether the point
$o\in F$ is singular or smooth. For the varieties of general
position, considered in this paper, all smooth points satisfy the
condition (f) whereas the singular points are of different
behavior. This point is discussed below in \S 3.

Assume first that $B=\mathop{\rm
centre}(\nu_E,V)\not\subset\mathop{\rm Sing}F$. In other words,
either the fiber $F$ is non-singular or $B$ is not a singular
point of this fiber. By the regularity condition, for any
irreducible subvariety $Y\subset F$ of codimension two (with
respect to $F$) the estimate (\ref{b9}) holds for a point $o\in B$
of general position. For all the multiplicities we have $\mu_i=1$.
The inequality (\ref{b9}) implies immediately the estimate
$$
\sum_{i\in J_s\cup J^+_m}p_im^h_i\leq 4n^2\Sigma_s.
$$
Since $m^h_i\leq m^h_1\leq 4n^2$, we get the inequality
\begin{equation}
\label{b13} \sum_{i\in J_l}p_im^h_i\leq 4n^2(\Sigma_s+\Sigma^-_m).
\end{equation}
This is the required estimate of singularities of the horizontal
component $Z^h$. Consider the vertical component $Z^v$. By the
condition (v) the inequality
\begin{equation}
\label{b10} m^v_i\leq m^v_1\leq \frac{2}{\mathop{\rm deg}V}d_v.
\end{equation}
holds. On the other hand, the generalized $K^2$-condition of depth
$\varepsilon$ implies the estimate
\begin{equation}
\label{b11} \frac{d_v}{\mathop{\rm
deg}V}<\frac{2en}{\nu_E(F)}+\varepsilon n^2.
\end{equation}
Combining (\ref{b10}) and (\ref{b11}), we obtain the inequality
$$
\sum_{i\in J_s\cup J^+_m}p_im^v_i<
2n\left(\frac{2e}{\nu_E(F)}+\varepsilon n
\right)(\Sigma_s+\Sigma^+_m).
$$
Taking into account that by definition
$\nu_E(F)=\sum\limits^k_{i=1}p_i\mu_i\geq\Sigma_s+\Sigma^+_m$, we
obtain finally
\begin{equation}\label{b12}
\sum_{i\in  J_s\cup J^+_m}p_i m^v_i< 4ne+2\varepsilon
n^2(\Sigma_s+\Sigma^+_m)
\end{equation}
Now the inequalities (\ref{a4}), (\ref{b13}) and (\ref{b12}) lead
to the following estimate:
$$
(4n^2(\Sigma_s+\Sigma^-_m)+4ne+2\varepsilon
n^2(\Sigma_s+\Sigma^+_m))(\Sigma_s+\Sigma_m+\Sigma_u)>
$$
$$
>((3\Sigma_s+2\Sigma_m+\Sigma_u)n+e)^2.
$$
Setting in this inequality $\varepsilon=2$ and taking into account
that $\Sigma_m=\Sigma^+_m+\Sigma^-_m$ (this is the crucial point),
after some easy arithmetic we get the inequality
$$
(n(\Sigma_s-\Sigma_u)+e)^2< 0.
$$
A contradiction.

The equality of the thresholds of canonical adjunction
$$
c_{\mathop{\rm virt}}(\Sigma)=c(\Sigma)
$$
is proved for the non-singular case $B\not\subset \mathop{\rm
Sing}F$.


\subsection{Estimating the multiplicities of a linear system: \\ the
singular case}

Now consider the case when $B=o\in F$ is a singular point of the
fiber. If the variety $F$ satisfies the condition (f) at the point
$o$, then the arguments of the previous section work well without
any modifications. If this is not the case, then $\mathop{\rm
deg}V=\mathop{\deg F}\geq 6$ (otherwise $\mathop{\rm deg}V\leq 4$
and the condition (f) holds automatically) and thus $\mathop{\rm
dim}F\geq 4$, $\mathop{\rm dim}V\geq 5$. Therefore the discrepancy
of the first exceptional divisor $E_1$ is at least 4 (we blow up a
smooth point $o\in V$). By assumption, the condition (fs) holds,
whence we get the inequality
$$
m^h_i\leq m^h_1\leq 3n^2,
$$
so that from the estimate (\ref{b2}) we obtain
$$
\sum_{i\in J_l}p_im^h_i\leq 3n^2(p_1+\Sigma_s+\Sigma^-_m).
$$
On the other hand, as we have mentioned above, $a(E_1,V)\geq 4$,
so that one can replace in the inequality (\ref{a4}) the numerator
by $(p_1+3\Sigma_s+2\Sigma_m+\Sigma_u)$. Now arguing as in Sec.
2.4 we get the inequality
$$
(n(\Sigma_s-\Sigma_u)+e)^2+(\Sigma_s-3p_1)(\Sigma_s+\Sigma_m+\Sigma_u)n^2
+
$$
$$
+ np_1(np_1+2(3\Sigma_s+2\Sigma_m+\Sigma_u)n+2e)<0.
$$
Obviously, $\Sigma_s-3p_1\geq-2p_1$ so that we again get a
contradiction.

The proof of the equality of the thresholds of canonical
adjunction $c_{\mathop{\rm virt}}(\Sigma)=c(\Sigma)$ is complete.

Because of the importance of this result for further work let us
formulate it as a separate claim.

{\bf Theorem 2.} {\it Assume that the fiber space $V/{\mathbb
P}^1$ satisfies the generalized $K^2$-condition of depth 2, the
conditions (v), (vs) and at least one of the conditions (f) or
(fs) at every point $o\in V$. Then for any movable linear system
$\Sigma\subset |-nK_V+lF|$ with $l\in{\mathbb Z}_+$ its virtual
and actual thresholds of canonical adjunction coincide:
$$
c_{\mathop{\rm virt}}(\Sigma)=c(\Sigma).
$$
}


\section{Varieties with a pencil of Fano double covers}

In order to apply the technique developed above to the fiber
spaces $V/{\mathbb P}^1$, the fibers of which are Fano double
hypersurfaces, one needs to describe movable linear systems on
these varieties and check the conditions on multiplicities of
horizontal and vertical cycles for all the families under
consideration: $1-8$ and $1^*-5^*$ of Sec. 0.2 and 0.3.

\subsection{Movable systems on the varieties of the type ((0),(2,0))}

Let $V/{\mathbb P}^1$ be a variety from the family ((0),(2,0)).
The construction of the variety $V$ can be alternatively described
in the following way.

Let $W_{\mathbb P}\subset{\mathbb P}$ be a hypersurface of degree
$2l$,
$$
\sigma_{\mathbb Y}\colon {\mathbb Y}\to{\mathbb P}
$$
the double cover branched over the divisor $W_{\mathbb P}$.
Consider the variety $Y={\mathbb P}^1\times{\mathbb Y}$, which is
realized as the double cover $\sigma_Y\colon Y\to X={\mathbb
P}^1\times {\mathbb P}$ branched over the divisor $W={\mathbb
P}^1\times W_{\mathbb P}$. Set $V=\sigma^{-1}_Y(Q)$, where
$Q\subset X$ is a smooth divisor of the type $(2,m)$. It is easy
to see that in this way we obtain the same variety $V$ as in Sec.
0.2.

By assumption, the divisor $Q\subset {\mathbb P}^1\times{\mathbb
P}$ is given by the equation
$$
A(x_*)u^2+2B(x_*)uv+C(x_*)v^2=0,
$$
where $A(\cdot),B(\cdot),C(\cdot)$ are homogeneous of degree $m$.
Here $(u:v)$ and $(x_*)=(x_0\colon\dots\colon x_{M+1})$ are
homogeneous coordinates on ${\mathbb P}^1$ and ${\mathbb P}$,
respectively.

Furthermore, let $H_{\mathbb P}$ be the class of a hyperplane in
${\mathbb P}$, $L_X=p^*_XH_{\mathbb P}$ the tautological class on
$X$, where $p_X\colon X\to{\mathbb P}$ is the projection onto the
second factor, $L_V=\sigma^*_YL_X|_V$. It is easy to see that
$$
K_V=-L_V,
$$
so that the anticanonical linear system $|-K_V|$ is free and
determines the projection $p_V=p_X\circ\sigma\colon V\to {\mathbb
P}$.

{\bf Lemma 3.1.} {\it The projection $p_V$ factors through the
double cover $\sigma_{\mathbb Y}\colon {\mathbb Y}\to {\mathbb
P}$. More precisely, there is a morphism $p\colon V\to {\mathbb
Y}$ such that
$$
p_V=\sigma_{\mathbb Y}\circ p.
$$
The degree of the morphism $p$ at a general point is equal to 2.}

{\bf Proof.} Consider a point $x\in{\mathbb P}\setminus W_{\mathbb
P}$ of general position. Set $\{y^+,y^-\}=\sigma^{-1}_{\mathbb
Y}(x)\subset{\mathbb Y}$. Set also
$$
L_x={\mathbb P}^1\times\{x\}\subset X, \quad L^{\pm}_x={\mathbb
P}^1 \times\{y^{\pm}\}\subset Y.
$$
It is obvious that the inverse image $\sigma^{-1}_Y(L_x)$ is the
disjoint union of the lines $L^+_x$ and $L^-_x$, whereas
$$
p_Y(L^{\pm}_x)=y^{\pm},
$$
where $p_Y\colon Y\to{\mathbb Y}$ is the projection onto the
second factor. The divisor $Q$ intersects $L_x$ at two distinct
(for a general point $x$) points $q_1,q_2$. Set
$$
\sigma^{-1}(q_i)=\{o^+_i,o^-_i \}\subset V,\quad o^{\pm}_i\in
L^{\pm}_x.
$$
The morphism $p$ is the restriction $p_Y|_V$. Obviously,
$$
p^{-1}(y^{\pm})=\{o^{\pm}_1,o^{\pm}_2\},
$$
where the sign $+$ or $-$ is the same in the right hand and left
hand side. This proves the lemma.

Let $\Delta\subset V$ be a subvariety of codimension 2, given by
the system of equations $A=B=C=0$. The subvariety $\Delta$ is
swept out by the lines $L_y={\mathbb P}^1\times\{y\}$ which are
contracted by the morphism $p$. Set $\Delta_{\mathbb
Y}=p(\Delta)$. Obviously,
$$
p\colon V \setminus \Delta \to {\mathbb Y} \setminus
\Delta_{\mathbb Y}
$$
is a finite morphism of degree 2. Let $\tau\in\mathop{\rm Bir}V$
be the corresponding Galois involution. It is easy to see that
$\tau$ commutes with the Galois involution $\alpha\in\mathop{\rm
Aut}V$ of the double cover $\sigma\colon V\to Q$, so that $\tau$
and $\alpha$ generate a group of four elements. Since the
involution $\tau$ is biregular outside the invariant closed subset
$\Delta$ of codimension 2, that is, $\tau\in\mathop{\rm
Aut}(V\setminus\Delta)$, the action of $\tau$ on the Picard group
$\mathop{\rm Pic}V$ is well defined.

Let $\Sigma\subset|-nK_V+lF|$ be a movable linear system.

{\bf Lemma 3.2.} (i) {\it The involution $\tau$ transforms the
pencil $|F|$ of fibers of the morphism $\pi$ into the pencil
$|mL_V-F|$.} (ii) {\it If $l<0$, then the involution $\tau$
transforms the linear system $\Sigma$ into the linear system
$$
\Sigma^+\subset|n^+L_V+l^+F|,
$$
where} $n^+=n+lm\geq 0,\,l^+=-l>0$.

{\bf Proof.} Obviously, $\tau^*L_V=L_V$. Let $F_t=\pi^{-1}(t)$ be
a fiber. We get
$$
p^{-1}(p(F_t))=F_t\cup\tau(F_t).
$$
However, $p(F_t)\sim mH_{\mathbb Y}=m\sigma^*_{\mathbb
Y}H_{\mathbb P}$ by the construction of the variety $V$. Since
$p^*H_{\mathbb Y}=L_V$, we obtain the claim (i). Thus
$\tau^*F=mL_V-F$. This directly implies the second claim of the
lemma.

{\bf Remark 3.1.} Varieties of the type $((0),(2,0))$ are similar
by their properties to the Fano fiber spaces considered in [19].

\subsection{Checking the $K$-condition}

Let us prove that the variety $V/{\mathbb P}^1$ from any of the
families $2-8$ or $1^*-5^*$ of Sec. 0.2, 0.3 satisfies the
$K$-condition. Fix a movable linear system $\Sigma\subset
|-nK_V+lF|$. We must show that $l\in{\mathbb Z}_+$. We will use
the arguments of the following two types:

\noindent (1) assume that there is a divisor  $E\subset V$, swept
out by a family of irreducible horizontal curves
$(C_{\delta},\delta\in \Delta)$ such that
$$
(-K_V\cdot C_{\delta})\leq 0.
$$
Since the linear system $\Sigma$ is movable we get $l\geq 0$,
which is what we need.

\noindent (2) Assume that there is a horizontal prime divisor $E$
on $V$ satisfying the inequality
$$
(-K_V\cdot E\cdot L^{M-1}_V)\leq 0.
$$
For a general divisor $D\in\Sigma$ the cycle $(D\circ E)$ is
effective and thus
$$
(D\cdot E\cdot L^{M-1}_V)\geq 0,
$$
whence, taking into account that $(H_F\cdot L^{M-1}_V)>0$, we get
$l\geq 0$.

Now let us consider varieties of the types $2-8$ and $1^*-5^*$ one
after another.

\subsubsection{Varieties of the type ((0),(1,1))}

Here $X={\mathbb P}^1\times {\mathbb P}$, $Q\subset X$ is given by
the equation $Au+ Bv=0$ where $A,B$ are homogeneous polynomials of
degree $m$ on ${\mathbb P}$. The pair of equations $A=B=0$ defines
on $Q$ the divisor $\Delta= {\mathbb P}^1\times \Delta_{\mathbb
P}$, which is swept out by horizontal lines, $\Delta_{\mathbb
P}\subset {\mathbb P}$ is the subvariety $\{A=B=0\}$ of
codimension two. For any line $L_x=\sigma^{-1}({\mathbb
P}^1\times\{x\}),x\in \Delta_{\mathbb P}$, we get
$$
(-K_V\cdot L_x)=(L_V\cdot L_x)=0.
$$
Now the remark (1) shows that the $K$-condition holds.

\subsubsection{Varieties of the type ((1),(0,1))}

Here
$$
{\cal E}={\cal O}^{\oplus(M+1)}_{{\mathbb P}^1}\oplus {\cal
O}_{{\mathbb P}^1}(1),
$$
and set $E\subset X$ to be the divisor of common zeros of all
sections $s\in {\cal O}_{{\mathbb P}^1}(1)$. Obviously,
$E={\mathbb P}^1\times {\mathbb P}^M$. Let $L_E=p^*_E(H_{\mathbb
P}|_{{\mathbb P}^M})$ be the tautological class on $E$, where
$p_E\colon E\to {\mathbb P}^M$ is the projection onto the second
factor. We get
$$
-K_V|_E=L_V|_E=L_E,
$$
$E\cap Q\sim mL_E$, so that $E\cap Q={\mathbb P}^1\times Q_E$,
where $Q_E\subset{\mathbb P}^M$ is a hypersurface of degree $m$.
Therefore, $V$ contains the divisor $\sigma^{-1}(E\cap Q)$, which
is swept out by the curves $L_x=\sigma^{-1}({\mathbb
P}\times\{x\}),x\in Q_E$. The anticanonical class $(-K_V)$ is
trivial on these curves. According to the remark (1), the
$K$-condition holds.

\subsubsection{Varieties of type ((2),(1,0))}

Here $-K_V=L_V-F$,
$$
{\cal E}={\cal O}^{\oplus(M+1)}_{{\mathbb P}^1}\oplus{\cal
O}_{{\mathbb P}^1}(2),
$$
set $E\subset X$ to be the divisor of common zeros of sections
$s\in {\cal O}_{{\mathbb P}^1}(2)$. Now
$$
\sigma^{-1}(E\cap Q)\in |L_V-2F|.
$$
It is easy to compute that
$$
(-K_V\cdot (L_V-2F)\cdot L^{M-1}_V)=2n(1-m)\leq 0.
$$
According to the remark (2), the $K$-condition holds.

\subsubsection{Varieties of the type ((2),(0,0))}

Here $-K_V=L_V$ and we can apply the remark (1): the divisor
$\sigma^{-1}(E\cap Q)\subset V$ (in the notations of the previous
case) is swept out by horizontal curves. The class $L_V$ is
trivial on these curves.

\subsubsection{Varieties of the type ((3),(0,0))}
Here
$$
{\cal E}={\cal O}^{\oplus(M+1)}_{{\mathbb P}^1}\oplus{\cal
O}_{{\mathbb P}^1}(3)
$$
and $-K_V=L_V-F$. Set $E\subset X$ to be the divisor of common
zeros of the sections $s\in{\cal O}_{{\mathbb P}^1}(3)$. Now the
divisor $\sigma^{-1}(E\cap Q)\subset V$ is swept out by horizontal
curves, on which $L_V$ is trivial. We apply the remark (1).

\subsubsection{Varieties of the type ((1,2),(0,0))}
Here $-K_V=L_V-F$ and there is a prime divisor $E\subset V$, such
that
$$
E\sim L_V-2F.
$$
It is easy to check that
$$
((L_V-2F)\cdot(L_V-F)\cdot L^{M-1}_V)=0.
$$
By the remark (2), the $K$-condition holds.

\subsubsection{Varieties of the type ((1,1,1),(0,0))}
Here $-K_V=L_V-F$ and moreover $\mathop{\rm codim}_V\mathop{\rm
Bs}|L_V-F|=3$. Thus for any pseudo-effective class $D\in
\mathop{\rm Pic}V$ we get the inequality
$$
(D\cdot(L_V-F)^2\cdot L^{M-2}_V)\geq 0.
$$
However, it is easy to compute that $((L_V-F)^3\cdot
L^{M-2}_V)=0$. Thus if $D=-nK_V+lF$ is just a pseudo-effective
class, then $l\in{\mathbb Z}_+$. In particular, the $K$-condition
holds.

The remaining five types of double spaces are considered in the
same way (with simplifications).

Q.E.D. for the $K$-condition for varieties of types $2-8$ and
$1^*-5^*$.

\subsection{Proof of birational rigidity}

Now in order to prove Theorem 1, we have to check the equality
$$
c_{\mathop{\rm virt}}(\Sigma)=c(\Sigma)
$$
for any movable linear system $\Sigma\subset|-nK_V+lF|$ with
$n\geq 1$ and $l\in{\mathbb Z}_+$. In its turn, by Theorem 2 it is
sufficient to verify that varieties of the types $1-8$ and
$1^*-5^*$ satisfy the generalized $K^2$-condition of depth 2 and
that a regular fiber space $V/{\mathbb P}^1$ satisfies the
conditions (v), (vs) and (f) or (fs) at every point (the condition
(h) follows directly from any of the two conditions (f) or (fs)).
The conditions (v) and (vs) are checked in [16], the condition (f)
for a regular smooth point of a fiber $o\in F$ is checked in [9].
It remains to show that at a double point $o$ the fiber $F$
satisfies at least one of the two conditions (f) or (fs).

If the singular point $o\in F$ lies outside the branch divisor of
the morphism $\sigma$, then the condition (f) is satisfied. This
is easy to check by means of the standard method of hypertangent
divisors. One should take into account that if $o\in F$ does not
lie on the branch divisor, then there are $l$ additional
hypertangent divisors, arising from the double cover. It is this
fact that makes it possible to obtain the estimate (\ref{b9}) for
any irreducible subvariety $Y\subset F$ of codimension 2 (with
respect to the fiber $F$). We do not give these arguments here
because they are standard and, in particular, parallel to another
estimate which is proved below for smooth points (where 4 is
replaced by 3).

Therefore let us assume that the singular point $o\in F$ lies on
the branch divisor. Let $Y\subset F$ be an irreducible subvariety
of codimension 2, $T=\sigma^{-1}(T_pG\cap G)$, $p=\sigma(o)$, the
tangent divisor. If $Y\not\subset T$, then considering the
effective cycle $(Y\circ T)$  of codimension 3 and applying the
standard technique of hypertangent divisors to this cycle (that
is, intersecting it with $D_i\in\Lambda_i$, $i=4,\dots,m-1$), we
obtain the inequality (\ref{b9}). Thus we may assume that
$Y\subset T$. Applying the technique of hypertangent divisors
(with $D_i\in \Lambda_i$ for $i=3,\dots,m-1$) we obtain
(\ref{b14}), the first of the two inequalities of the condition
(fs). Let us prove the second one.

Let $\varphi\colon \widetilde T\to T$ be the blow up of the point
$o$, $E=\varphi^{-1}(o)\subset \widetilde T$ the exceptional
divisor. It is easy to see that the double cover $\sigma$ presents
$E$ as the double cover of the quadric $E_G$ branched over the
divisor $\widetilde W\cap E_G$, that is, the section of $E_G$ by a
quadric hypersurface. The quadric $E_G$ is the exceptional divisor
of the blow up of the point $p=\sigma(o)$ on the variety $T_pG\cap
G$, $\widetilde W$ denotes the strict transform of the restriction
$W|_G$. By assumption, $Y$ is a prime divisor on $T$. Therefore,
$Y_E=(\widetilde Y\circ E)$ is an effective divisor on the double
quadric $E$, and $\mathop{\rm deg}Y_E=\mathop{\rm
mult}\nolimits_oY$. Let $x\in E$ be an arbitrary point.

Set $q=\sigma(x)\in E_G$, $R=\sigma^{-1}(T_qE_G\cap E_G)$. The
divisor $R\subset E$ is irreducible, $\mathop{\rm
deg}R=4,\mathop{\rm mult}\nolimits_xR=2$. If $Z\subset E$ is a
prime divisor, different from $R$, then the cycle $(Z\circ R)$ of
codimension 2 on $E$ is well defined, and moreover
$$
\mathop{\rm deg}Z=\mathop{\rm deg}(Z\circ R)\geq\mathop{\rm
mult}\nolimits_x(Z\circ R)\geq 2\mathop{\rm mult}\nolimits_xZ.
$$
This implies that
$$
\mathop{\rm mult}\nolimits_x\widetilde Y\leq \mathop{\rm
mult}\nolimits_xY_E\leq \frac12\mathop{\rm
deg}Y_E=\frac12\mathop{\rm mult}\nolimits_oY.
$$
This proves the second inequality of the condition (fs). As for
the generalized $K^2$-condition of depth 2, it follows from the
inequality
$$
((K^2_V-2H_F)\cdot L^{M-1}_V)\leq 0,
$$
which is easy to check for varieties of types $1-8$ and  $1^*-5^*$
of Sec. 0.2 and 0.3.

Q.E.D. for Theorem 1.


\subsection{Multiplicities of subvarieties of codimension 2}

In the remaining part of this section we prove the following
claim.

{\bf Proposition 3.1.} {\it Let $o\in F$ be a smooth point of the
fiber lying outside the branch divisor of the morphism $\sigma$.
Then for any irreducible subvariety $Y\subset F$ of codimension
two (and thus for any effective cycle of pure codimension two) the
following estimate is true:}

\begin{equation}
\label{c2} \frac{\mathop{\rm mult}\nolimits_o}{\mathop{\rm
deg}}Y\leq \frac{3}{\mathop{\rm deg}V}=\frac{3}{2m}.
\end{equation}

{\bf Remark 3.2.} Proposition 3.1 outlines another approach to
proving birational rigidity of the varieties under consideration,
based on Proposition 1.7. If a smooth point $o\in F$ lies on the
branch divisor of the morphism $\sigma$, then there exists an
irreducible subvariety $Y^+\subset F$ of codimension 2 for which
the estimate of the condition (f) cannot be improved: the linear
system $|H-2o|$ is movable in this case and spanned by two
divisors, $\sigma^*(T_pG\cap G)$ and $\sigma^*(T_pW\cap G)$, where
$p=\sigma(o)$. The intersection of these divisors gives a
subvariety $Y^+$ with this property. It can be shown, however,
that for any other subvariety of codimension two, $Y\neq Y^+$ the
inequality (\ref{c2}) holds. As for the subvariety $Y^+$, it
satisfies the estimate (\ref{c2}) at all other points, including
the infinitely near points of the first order lying over the point
$o$. Now using the discrepancy arguments one can prove birational
rigidity using only Proposition 1.7. However, a combination of
these ideas with Theorem 2 makes it possible to improve the
technique even further. We will use this strategy in the
subsequent papers.

{\bf Proof of Proposition 3.1.} Arguing in the same way as in
[16], \S 2, we give two methods of obtaining the estimate
(\ref{c2}). The first one is much more simple, however works only
for sufficiently high dimensions $M\geq M_0$. The second method
gives the estimate (\ref{c2}) for all dimensions, however it
requires additional conditions of general position which one has
to substantiate (Sec. 3.5). Whichever method is used, the first
step is the same.

{\bf Lemma 3.3.} {\it Assume that the inequality
\begin{equation} \label{c6}
\frac{\mathop{\rm mult}\nolimits_o}{\mathop{\rm
deg}}Y>\frac{3}{2m}
\end{equation}
holds. Then $\sigma(Y)\subset T_pG$, where} $p=\sigma(o)$.

{\bf Proof.} Assume the converse: $\sigma(Y)\not\subset T_pG$.
Then the intersection $Y\cap T$, where $T=\sigma^*(T_pG\cap G)$ is
the tangent section of the fiber $F$ at the point $o$, is of
codimension 3 and therefore the effective cycle $(Y\circ T)$ is
well defined. For some component $Y_3$ of this cycle the
inequality
\begin{equation}\label{c3}
\frac{\mathop{\rm mult}\nolimits_o}{\mathop{\rm
deg}}Y_3>\frac{3}{m}
\end{equation}
holds. Now arguing in the usual way let us consider general
divisors of the hypertangent linear systems
$$
D_i\in \Lambda_i,\quad  D^+_j\in \Lambda_j,
$$
where
$$
i\in\{4,\dots,m-1\},\quad j\in\{l,\dots,2l-2\}
$$
for $m\leq 2l$ and
$$
i\in\{4,\dots,m-2\},\quad j\in\{l,\dots,2l-1\}
$$
for $m\geq 2l+1$. Taking into account the codimension of the base
set of the system $\Lambda_i$, we see that the set
$$
Y_3\cap(\bigcap\limits_iD_i)\cap(\bigcap\limits_jD^+_j)
$$
is one-dimensional in a neighborhood of the point $o$. From this
we obtain in the standard way the estimate
\begin{equation}\label{c4}
\frac{\mathop{\rm mult}\nolimits_o}{\mathop{\rm deg}}Y_3\leq
\frac{4l}{m(2l-1)}
\end{equation}
for $m\leq 2l$ and the estimate
\begin{equation}\label{c5}
\frac{\mathop{\rm mult}\nolimits_o}{\mathop{\rm
deg}}Y_3\leq\frac{2}{m-1}
\end{equation}
for $m\geq 2l+1$. In any case the estimates (\ref{c4}) and
(\ref{c5}) are incompatible with the inequality (\ref{c3}). Q.E.D.
for the lemma.

Now let us prove Proposition 3.1. If the estimate (\ref{c2}) does
not hold, that is, if the estimate (\ref{c6}) holds, then by what
was proved above, $\sigma(Y)\subset T_pG$, that is,
$$
Y\subset T_1=\sigma^{-1}(G\cap T_pG).
$$
By the regularity conditions, $T_1\subset F$ is an irreducible
subvariety with the isolated double point $o\in T_1$. Consider the
closed set
$$
T_{12}=\sigma^{-1}(\sigma(T_1)\cap T_p\sigma(T_1)).
$$
By the regularity conditions, $T_{12}\subset F$ is an irreducible
subvariety of codimension two, and moreover $\mathop{\rm
mult}\nolimits_oT_{12}=6$, so that the equality
$$
\frac{\mathop{\rm mult}\nolimits_o}{\mathop{\rm
deg}}T_{12}=\frac{3}{2m}
$$
holds, whence by the inequality (\ref{c6}) we get, that $Y\neq
T_{12}$. As above, let
$$
D_i\in \Lambda_i,\,D^+_j\in\Lambda^+_j
$$
be general hypertangent divisors, where
$$
i\in I=\{2,4,\dots,m-1\},\quad j\in J=\{l,\dots,2l-2\}
$$
(in $I$ we omit the element $i=3$) for $m\leq 2l$ and
$$
i\in I\{2,4,\dots,m-2\},\quad j\in J=\{l,\dots,2l-1\}
$$
(in $I$ we again omit $i=3$) for $m\geq 2l+1$. By the regularity
conditions the intersection
$$
Y\cap(\bigcap\limits_{i\in I}D_i)\cap(\bigcap\limits_{j\in
J}D^+_j)
$$
is one-dimensional in a neighborhood of the point $o$, whence we
get the estimates
$$
\frac{\mathop{\rm mult}\nolimits_o}{\mathop{\rm
deg}}Y\leq\frac{8l}{3m(2l-1)}
$$
and
$$
\frac{\mathop{\rm mult}\nolimits_o}{\mathop{\rm
deg}}Y\leq\frac{4}{3(m-1)}
$$
for $m\leq 2l$ and $m\geq 2l+1$, respectively. It is easy to check
that these estimates prove Proposition 3.1 in all cases except for
the following ones:
$$
m\in\{3,4,5,6,7,8,\},\quad l\in \{3,4\}.
$$

Now let us give a more refined way of arguing, which requires
strong regularity conditions but works in all dimensions including
low ones. This method is completely similar to the method by means
of which the condition (vs) was proved in [16] for a double point
of the fiber $o\in F$ outside the branch divisor. For this reason
here we just describe the main steps of the proof, emphasizing the
changes which should be made in the arguments of the paper [16].

So let us assume that for an irreducible variety $Y\subset F$ the
inequality (\ref{c6}) holds, which contradicts the required
estimate (\ref{c2}). By Lemma 3.3, $Y\subset T_1$. The point $o$
is an isolated factorial singularity of the variety $T_1$.
Moreover, $\mathop{\rm Pic}T_1={\mathbb Z}H_T$, where
$$
H_T=H_F|_{T_1}
$$
is the class of a hyperplane section. Let
$\varphi_T\colon\widetilde T\to T_1$ be the blow up of the point
$o$, $E_T\subset \widetilde T$ the exceptional divisor (it is
irreducible), $\widetilde Y\subset \widetilde T$ the strict
transform of the divisor $Y$. We get
$$
\widetilde Y\sim \alpha H_T-\beta E_T,
$$
where $\beta/\alpha>3/2$ by the estimate (\ref{c6}).

{\bf Lemma 3.4.} {\it The prime divisor $T_{12}\subset T_1$ is
swept out by a family of curves $\{C_{\delta},\delta\in \Delta\}$,
the general curve of which is irreducible and satisfies the
estimate}
$$
\frac{\mathop{\rm mult}\nolimits_o}{\mathop{\rm
deg}}C_{\delta}>\frac23.
$$
Assume that the lemma is proved and consider the strict transform
$\widetilde C_{\delta}\subset \widetilde T$ of the general curve
$C_{\delta},\delta\in \Delta$. Obviously,
$$
(\widetilde Y\cdot \widetilde C_{\delta})< 0,
$$
so that $C_{\delta}\subset Y$ and thus $T_{12}\subset Y$.
Consequently, $Y=T_{12}$ which is impossible. Q.E.D. for
Proposition 3.1.

To prove Lemma 3.4, one needs the arguments which are absolutely
similar to the arguments that were used in the proof of Lemma 2.2
in [16], so that we will not repeat them, just remind the main
steps. For $m\geq 4$, $l\geq 3$ define the sequence of integers
$$
c_e= \sharp [4,e]\cap{\cal M}+ \sharp [3,e]\cap {\cal L},\quad
e\in {\mathbb Z}_+,
$$
where ${\cal M}=\{2,\dots,m-1\}$, ${\cal L}=\{l,\dots,2l-1\}$, and
construct the {\it ordering function}
$$
\chi\colon \{1,\dots,m+l-4\}\to {\mathbb Z}_+,
$$
setting $\chi([c_{e-1}+1,c_e]\cap {\mathbb Z}_+)=e$. By the
regularity condition, for a general set of hypertangent divisors
$$
{\mathbb D}=\{D_i\in \Lambda^P_{\chi(i)},\,\,\,i=1,\dots,m+l-4\}
\in \Lambda^P
$$
the closed algebraic set
$$
R_i({\mathbb D})=\mathop{\bigcap}\limits^i_{j=1} D_i\cap T_P
$$
is of codimension $i$ in $T_{12}$ for $i=1,\dots,m+l-4$. Thus we
get an effective 1-cycle
$$
R({\mathbb D})=R({\mathbb D})=(T_P\circ D_1 \circ \dots \circ
D_{m+l-4})= \sum_{\delta_i\in\Delta}C_{\delta_i}+\Phi ,
$$
where $(C_{\delta},\delta \in\Delta)$ is a movable family of
curves, sweeping out $T_{12},\Phi$ is the fixed part of the family
$R({\mathbb D})$. It is easy to see that $\Phi$ is exactly the
1-cycle of lines on $F$, passing through the point $o\in F$, that
is, $\mathop{\rm deg}\Phi=\mathop{\rm deg}{\mathbb L}(o)$.
Therefore,
$$
\frac{\mathop{\rm mult}\nolimits_o}{\mathop{\rm
deg}}C_{\delta}\geq \frac{\displaystyle\frac{m!}{4}\cdot
\frac{(2l)!}{l!}-\lambda_{m,l} }{\displaystyle\frac{2m!}{3}\cdot
\frac{(2l-1)!}{(l-1)!}-\lambda_{m,l} }>\frac23
$$
by Proposition 3.2 (which is proved below). In a similar way we
argue when $m=3$ or $l=2$. Proof of Lemma 3.3 is complete.

Let us emphasize that in contrast to the situation considered in
[16], here the hardest point is to estimate the number of lines
(taken with multiplicities), passing through the point $o\in F$.
In [16] the double points of the fibers are considered. There are
but finitely many such points so that for a general Fano fiber
space $V/{\mathbb P}^1$ all the lines passing through a singular
point of a fiber are of multiplicity one, and this task becomes
trivial. However, in out case $o\in F$ is an arbitrary point of a
fiber. As we will show just now, the multiplicity of a line
passing through some specially chosen point $o\in F$ can be very
high. But not high enough to prevent the method of proving Lemma
2.2 in [16] to work in our case. All the other differences between
the proof of Proposition 3.1 of this paper and that of Proposition
2.3 in [16] are inessential.


\subsection{Estimating the number of lines}

For a point $o\in V\setminus\sigma^{-1}(W)$ outside the branch
divisor we define the {\it algebraic cycle of lines} ${\mathbb
L}(o)$ on $V$ passing through $o$ as the $0$-cycle of the
subscheme
$$
\{q_1=\dots=q_m=g_{l+1}=\dots=g_{2l}=0\}
$$
on ${\mathbb E}={\mathbb P}(T_p{\mathbb P})\cong{\mathbb P}^M.$

{\bf Proposition 3.2.} {\it For a general (in the sense of Zariski
topology) fiber space $V/{\mathbb P}^1$ for any smooth point $o\in
V\setminus\sigma^{-1}(W)$ the following estimate holds:}
$$
\mathop{\rm deg}{\mathbb L}(o)\leq
\lambda_{m,l}=\frac{m!}{6}\frac{(2l-1)!}{(l-1)!}-1
$$

{\bf Proof.} We will describe the scheme of arguments and give
with all details the main technical lemma which makes it possible
to estimate the multiplicities of lines passing through an
arbitrary point $o\in V$. The very computations, which are, on one
hand, elementary, and on the other hand, tedious, will not be
given.

Set
$$
{\cal H}=\prod^m_{i=1}{\mathbb P}(H^0({\mathbb P}^m,{\cal
O}_{{\mathbb P}^M}(i)))\times\prod^{2l}_{j=l+1}{\mathbb
P}(H^0({\mathbb P}^M,{\cal O}_{{\mathbb P}^M}(j))).
$$
To simplify the notations, we write down a set of non-zero
polynomials $(q_{\sharp},g_{\sharp})\in{\cal H}$ as
$$
h_{\sharp}=(h_1,\dots,h_{M+1})\in{\cal H},
$$
where interchanging the direct factors we assume that $\mathop{\rm
deg}h_{i+1}\geq\mathop{\rm deg}h_i$. Let
$$
{\cal H}^+\subset  {\cal H}
$$
be the space of all collections $h_{\sharp}$ such that the scheme
of common zeros of the polynomials $h_i$ is zero-dimensional or
empty. For $h_{\sharp}\in{\cal H}^+$ let ${\mathbb L}(h_{\sharp})$
be the $0$-cycle of common zeros of the system $h_{\sharp}$. Now
we have the following fact.

{\bf Proposition 3.3.} {\it The codimension of the closed set
$$
{\cal H}^+(j)=\{h_{\sharp}\in{\cal H}^+|\mathop{\rm deg}{\mathbb
L}(h_{\sharp})\geq j\}
$$
for $j=\lambda_{m,l}+1$ is not less than} $M+2$.

Proposition 3.2 follows from this fact automatically. Now we
explain how to prove Proposition 3.3. In fact, for
$j=\lambda_{m,l}+1$ a much sharper bound for the codimension
$\mathop{\rm codim}{\cal H}^+(j)$ can be obtained in this way than
$M+2$ but we do not need that.

First of all, note an obvious fact: for $j=0,1,\dots,M+2$ the
following estimate holds:
$$
\mathop{\rm codim}\{h_{\sharp}\in{\cal H}^+|\,{\sharp}\mathop{\rm
Supp}{\mathbb L}(h_{\sharp})\geq j\}=j.
$$
This implies, that we can assume that the support of the cycle
${\mathbb L}(h_{\sharp})$ consists of $j\leq M+1$ distinct points
(that is, at most $M+1$ distinct lines pass through any point
$o\in V$ in the fiber $F\ni o$). However, as the dimension $M$
grow, the multiplicity of these points can be very high. Let us
consider the following

{\bf Example.} Let $p\in S=\{h_1=\dots=h_{M-a}=0\}$ be a
non-singular point, so that $S\ni p$ is a germ of a smooth
$a$-dimensional variety. When we require that
\begin{equation}\label{c7}
\mathop{\rm mult}\nolimits_p\{h_j|_S=0\}\geq 2,
\end{equation}
$j=M-a+1,\dots,M+1$, we impose precisely $a+1$ conditions on
$h_j$. Taking into account that the point $p$ is arbitrary, we
obtain
$$
a(a+1)\sim a^2
$$
conditions on the polynomials $h_j$, $M-a+1\leq j\leq M+1$, and
these conditions are independent. Thus for $a\sim \sqrt{M}$ there
are points on $V$ which satisfy the condition (\ref{c7}). For
those points we get
$$
\mathop{\rm deg}{\mathbb L}(o)\geq 2^a\sim 2^{\sqrt{M}}.
$$
This example is a model one showing how the function
$$
\mathop{\rm max}\limits_{o\in V}\mathop{\rm deg}{\mathbb L}(o)
$$
grows for a general variety $V/{\mathbb P}^1$ as $M\to \infty$.
Note, however, that for small values of $M$ we obtain the same
estimate as for the number of distinct lines
$$
\mathop{\rm max}_{o\in V}\,{\sharp}\mathop{\rm Supp}{\mathbb
L}(o).
$$

Now let us consider the problem of estimating the codimension of
the closed set ${\cal H}^+(j)$. By what was said above, we may
make the set ${\cal H}^+$ smaller and assume that
$$
{\sharp}\mathop{\rm Supp}{\mathbb L}(h_{\sharp})\leq M+1
$$
for any $h_{\sharp}\in{\cal H}^+$. For each $i=1,\dots,M$ set
$$
Y_i=\{h_1=\dots=h_i=0\}=\mathop{\circ}\limits^i_{j=1}\{h_j=0\}
$$
to be the algebraic cycle of the scheme-theoretic intersection of
the hypersurfaces $\{h_j=0\}$, $\mathop{\rm codim}_{\mathbb
E}Y_i=i$. Set also
$$
Y_{M+1}=\{h_1=\dots=h_{M-1}=h_{M+1}=0\}.
$$
Fix the points $o_1,\dots,o_k\in{\mathbb E}$, $k\leq M+1$ and for
each integer-valued matrix
$$
{\mathbb M}=\|\mu_{ij}\|_{\scriptstyle \begin{array}{l}
\scriptstyle 1\leq i\leq k,\\ \scriptstyle 1\leq j\leq M+1
\end {array}}\in \mathop{\rm Mat}\nolimits_{k\times(M+1)}({\mathbb Z}_+)
$$
consider the set
$$
{\cal H}({\mathbb M})=\{h_{\sharp}|\mathop{\rm
mult}\nolimits_{o_i}Y_j=\mu_{ij}\}\subset{\cal H}^+.
$$
(Naturally, it is sufficient to consider only those matrices which
satisfy the non-decreasing condition
$$
\mu_{i,1}\leq \mu_{i,2}\leq \dots \leq \mu_{i,M-1}\leq \mathop{\rm
min}(\mu_{i,M},\mu_{i,M+1}),
$$
otherwise either ${\cal H}({\mathbb M})$ is empty or the point
$o_i$ can be removed.) Obviously,
$$
\mathop{\rm deg}{\mathbb L}(h_{\sharp})\leq
\sum^k_{i=1}\mathop{\rm min}(\mu_{i,M},\mu_{i,M+1}).
$$
Now to prove Proposition 3.3 we must estimate the codimension
$$
\mathop{\rm codim}\nolimits_{\cal H} \overline{{\cal H}({\mathbb
M})}
$$
for all matrices ${\mathbb M}$, satisfying the inequality
$$
\sum^k_{i=1}\mathop{\rm min}(\mu_{i,M},\mu_{i,M+1})\geq
\lambda_{m,l}+1.
$$
Below in Sec. 3.6 we give an inductive method of estimating the
codimension $\mathop{\rm codim}\overline{{\cal H}(M)}$ at each
step of our procedure of making $Y_{i+1}$ from $Y_i$. Now the
proof of Proposition 3.3 is completed by tedious but absolutely
elementary computations based on Lemma 3.5, which is proved below.
We do not give these computations here.

Q.E.D. for Propositions 3.3 and 3.2.


\subsection{A method of estimating the degree}

Let $Y\subset {\mathbb A}={\mathbb C}^N$ be an effective cycle of
dimension $a\geq 1$, $o_1,\dots,o_k\in Y$ a set of pair-wise
distinct points. Set
$$
\mu_i=\mathop{\rm mult}\nolimits_{o_i}Y,
$$
$i=1,\dots,k$. Let ${\cal P}_e(z_1,\dots,z_N)$ be the space of all
(non-homogeneous) polynomials of degree not higher than $e\in
{\mathbb Z}_+,$ in the variables $z_1,\dots,z_N$. Let
$c_1,\dots,c_k\in{\mathbb Z}_+$ be a set of positive integers such
that
$$
c_1+\dots+c_{k-1}+(k-1)\leq e,\quad (k-1)a\leq N
$$
(if $k=1$, then no restrictions are imposed). Set
$$
a_1=\dots=a_{k-1}=a, \quad a_k=\mathop{\rm min}(a,N-(k-1)a),
$$
$c^*_i=c_i$ for $i=1,\dots,k-1$,
$$
c^*_k=\mathop{\rm min}(c_k,e-c_1-\dots-c_{k-1}-(k-1)).
$$
Let
$$
U_Y=\{f\in{\cal P}_e(z_*)|\mathop{\rm dim} \{f=0\}\cap\mathop{\rm
Supp}Y=a-1\}
$$
be an (open) set of polynomials that do not vanish identically on
each component of the cycle $Y$. Set
$$
U_{c_*}=U{(c_1,\dots,c_k)}=\{f\in U_Y|\mathop{\rm
mult}\nolimits_{o_i}\{f|_Y=0\}\geq c_i\mu_i+1\}.
$$
Obviously, the subset $U_{c_*}\subset U_Y$ is a Zariski closed
set.

{\bf Lemma 3.5.} {\it The following estimate holds:}
\begin{equation} \label{c8}
\mathop{\rm codim}U_{c_*}\geq
\Delta(c_1,\dots,c_k)=\sum^k_{i=1}{{a_i+c^*_i}\choose c^*_i}.
\end{equation}

{\bf Proof.} Let $\overline{U}_{c_*}\subset {\cal P}_e(z_*)$ be
the closure of the set $U_{c_*}$. In order to prove the inequality
(\ref{c8}), it is sufficient to produce a closed subset $Z\subset
{\cal P}_e(z_*)$, satisfying the inequality
\begin{equation} \label{c9}
\mathop{\rm dim}Z=\Delta(c_1,\dots,c_k),
\end{equation}
and such that
$$
Z\cap \overline U_{c_*}=\{0\}.
$$
Obviously, $0\in \overline U_{c_*}$, since $\overline U_{c_*}$ is
a cone with the vertex at zero. We first explain how such a set
can be constructed in the case when there is only one point
$o=o_1$, $k=1$, which we without loss of generality can assume to
be the origin $o=(0,\dots,0)\in {\mathbb A}$. To simplify our
notations, we write
$$
c_1=c^*_1=c\leq e,\quad a_1=a\leq N, \quad \mu_1=\mu.
$$
Let $\varphi\colon\widetilde{\mathbb A}\to {\mathbb A}$ be the
blow up of the point $o$, $E=\varphi^{-1}(o)\subset
\widetilde{\mathbb A}$ the exceptional divisor $E\cong{\mathbb
P}^{N-1}$. Consider an effective divisor $D\subset{\mathbb A}$,
$o\in D$, which contains no component of the cycle $Y$, so that
$$
\mathop{\rm dim}(\mathop{\rm Supp}Y\cap\mathop{\rm Supp}D)=a-1
$$
and the effective cycle
$$
Y_D=(Y\circ D)
$$
of dimension $a-1$ is well defined. Recall how the multiplicity of
this cycle at the point $o$ is computed. Let $\widetilde Y$,
$\widetilde D\subset \widetilde{\mathbb A}$ be the strict
transforms of the cycle $Y$ and divisor $D$, respectively. Set
$$
Y_E=(\widetilde Y\circ E)= \sum_{i\in I}b_iB_i
$$
to be the projectivized tangent cone to $Y$ at the point $o\in
\mathop{\rm Supp}Y$. We get
$$
(\widetilde D\circ \widetilde Y)=(\widetilde{D\circ Y})+\sum_{i\in
I}d_iB_i
$$
for some nonnegative $d_i\in{\mathbb Z}_+$. Now
$$
\mathop{\rm mult}\nolimits_ o Y_D=\mu\mathop{\rm mult}\nolimits_ o
D+\sum_{i\in I}d_i\mathop{\rm deg}B_i.
$$
(It is a standard fact of the intersection theory, see [21].) It
follows that if $B_i\not\subset\widetilde D$ for all $i\in I$,
then $d_i=0$ and thus the following equality holds:
$$
\mathop{\rm mult}\nolimits_oY_D=\mu\mathop{\rm mult}\nolimits_o D.
$$
We say that a set of linear functions $z_1,\dots,z_a$ is {\it
correct} with respect to the pair $(o\in Y)$, if

\begin{itemize}

\item $B_i\not\subset\widetilde{(z_1)}$ for all $i\in I$;

\item the projection $\pi_{(z_1,\dots,z_a)}\colon {\mathbb C}^N\to{\mathbb
C}^a$ is dominant on each irreducible component of the cycle $Y$;

\item the following infinitesimal condition is satisfied at the
point $o$. Taking $(z_1\colon\dots\colon z_N)$ for homogeneous
coordinates on $E$, consider the affine set $\{z_1\neq 0\}$ with
the coordinates $(y_2,\dots,y_N)$, where $y_i=z_i/z_1$. The
collection of functions $(y_2,\dots,y_a)$ determines a projection
$\iota\colon{\mathbb C}^{N-1}\to {\mathbb C}^{a-1}$. Now our
condition is formulated in the following way: the restriction
$$
\iota|_{B_i}\colon B_i\cap\{z_1\neq o\}\to{\mathbb C}^{a-1}
$$
of the projection $\iota$ onto each component $B_i$ is a dominant
map.
\end{itemize}
It is easy to see that a general set $(z_1,\dots,z_a)$ of linear
functions is correct with respect to the pair $(o\in Y)$. By the
definition of correctness for any non-zero polynomial $f\in{\cal
P}_e(z_i,\dots,z_a)$ in the variables $z_1,\dots,z_a$ we get
$$
B_i\not\subset \widetilde{(f)}
$$
and thus by what was said above $\mathop{\rm
mult}\nolimits_oY_{(f)}\leq\mu\mathop{\rm deg}f$. Thus we can set
$$
Z={\cal P}_c(z_1,\dots,z_a).
$$
For any non-zero polynomial $f\in Z$ the set $\{f=0\}$ contains
entirely none of the irreducible components of the cycle $Y$, so
that $Z\setminus\{0\}\subset U_Y$. Therefore $Z\cap\overline
U_c=\{0\}$, as required.

Now let us consider the general case of an arbitrary $k\geq 2$.
Here we have $k$ points $o_1,\dots,o_k\in\mathop{\rm Supp}Y$.
Consider a system of affine functions
$$
\begin{array}{ccc}
l_{1,1},&\dots,& l_{1,\,a_1},\\
&\dots&\\
l_{k,1},&\dots,& l_{k,\,a_k},
\end{array}
$$
$a_1=\dots=a_{k-1}=a$, $a_k=\mathop{\rm min}(a,N-(k-1)a)$,
satisfying the following conditions:

\begin{itemize}
\item
the linear parts of the affine functions $l_{*,*}$ are linear
independent, that is, form a part of a basis of the space
${\mathbb C}^N$;

\item for any $i\in\{1,\dots,k-1\}$ we have
$l_{i,\alpha}(o_i)=0$ and the system of functions
$(l_{i,1},\dots,l_{i,a})$ is correct with respect to the pair
$(o_i\in Y)$; $l_{k,\alpha}(o_k)=0$ and the system of functions
$(l_{k,1},\dots,l_{k,a_k})$ is a part of a correct set for the
pair $(o_k\in Y)$;

\item for $i\neq j$ we have $l_{i,1}(o_j)\neq 0$.

\end{itemize}

Now set
\begin{equation}\label{c10}
Z=\sum^k_{i=1}\left(\prod^{i-1}_{j=1}l^{c_j+1}_{j,1}\right){\cal
P}_{c^*_i}(l_{i,1},\dots,l_{i,a_i}),
\end{equation}
where ${\cal P}_\alpha(\sharp)$ denotes the linear space of all
polynomial functions of degree $\alpha$ in the affine functions
$\sharp$. Note that the sum in (\ref{c10}) is direct, since the
linear parts of the functions $l_{\sharp}$ are linear independent.

{\bf Lemma 3.6.} {\it The closed algebraic subsets $Z$ and
$\overline{U}_{c_*}$ intersect each other by zero only.}

{\bf Proof.} For $f\in Z$ we have the decomposition
$$
f=f_1+\dots+f_k,\quad f_i\in
Z_i=\left(\prod\limits^{i-1}\limits_{j=1}l^{c_j+1}_{j,1}\right){\cal
P}_{c^*_i}(l_{i,1},\dots,l_{i,a_i}),
$$
which is uniquely determined since
$Z=\mathop{\oplus}\limits^k_{i=1}Z_i$.

By construction, for each $i\in\{1,\dots,k\}$ we have
\begin{equation} \label{c11}
\mathop{\rm mult}\nolimits_{o_j}\{f_i|_Y=0\}\geq (c_j+1)\mu_j\geq
c_j\mu_j+1
\end{equation}
for all $j\leq i-1$. In particular, if $f\in
Z\cap\overline{U}_{c_*}$, then the following estimate holds:
$$
\mathop{\rm mult}\nolimits_{o_1}\{f_1|_Y=0\}\geq c_1\mu_1+1,
$$
since for $f_2,\dots,f_k$ the estimate (\ref{c11}) with $j=1$ is
satisfied. Arguing as above, we conclude that $f_1\equiv 0$.
Assume that it is already proved that
$$
f_1\equiv\dots\equiv f_{\gamma}\equiv 0.
$$
Since for $j=\gamma+1$ the estimate (\ref{c11}) holds for all
$i\geq \gamma+2$, we conclude that
\begin{equation}\label{c12}
\mathop{\rm mult}\nolimits_{o_{\gamma+1}}\{f_{\gamma+1}|_Y=0\}\geq
c_{\gamma+1}\mu_{\gamma+1}+1.
\end{equation}
However by construction
$$
f_{\gamma+1}=\left(\prod^{\gamma}_{j=1}l^{c_j+1}_{j,1}\right)f^{\sharp}_{\gamma+1},
$$
where $f^{\sharp}_{\gamma+1}\in{\cal
P}_{c^*_{\gamma+1}}(l_{\gamma+1,1},\dots,l_{\gamma+1,a_{\gamma+1}})$.
Taking into account that for $j=1,\dots,\gamma$ we have
$l_{j,1}(o_{\gamma+1})\neq 0$, we see that the estimate
(\ref{c12}) remains valid if $f_{\gamma+1}$ is replaced by
$f^{\sharp}_{\gamma+1}$. Now arguing like in the case $k=1$ above,
we conclude that $f_{\gamma+1}\equiv 0$. Q.E.D. for Lemma 3.6.

Taking into account that the equality (\ref{c9}) holds (this is
obvious from the explicit construction of the linear space $Z$
(\ref{c10})), we complete the proof of Lemma 3.5.

\section*{References}
{\small

\noindent 1. Brown G., Corti A. and Zucconi F. Birational geometry
of 3-fold Mori fibre spaces. Preprint, 2003, 40 p. arXiv:
math.AG/0307301. \vspace{0.3cm}

\noindent 2. Corti A., Pukhlikov A. and Reid M., Fano 3-fold
hypersurfaces, in ``Explicit Birational Geometry of Threefolds'',
London Mathematical Society Lecture Note Series {\bf 281} (2000),
Cambridge University Press, 175-258. \vspace{0.3cm}

\noindent 3. Grinenko M.M., Birational automorphisms of a
three-dimensional double cone. Sbornik: Mathematics. {\bf 189}
(1998), no. 7, 37-52. \vspace{0.3cm}

\noindent 4. Grinenko M.M., Birational properties of pencils of
del Pezzo surfaces of degrees 1 and 2. Sbornik: Mathematics. {\bf
191} (2000), no. 5, 17-38. \vspace{0.3cm}

\noindent 5. Grinenko M.M., Birational properties of pencils of
del Pezzo surfaces of degrees 1 and 2. II. Sbornik: Mathematics.
{\bf 194} (2003). \vspace{0.3cm}

\noindent 6. Iskovskikh V.A. and Manin Yu.I., Three-dimensional
quartics and counterexamples to the L\" uroth problem, Math. USSR
Sb. {\bf 86} (1971), no. 1, 140-166. \vspace{0.3cm}

\noindent 7. Pukhlikov A.V., Birational automorphisms of
three-dimensional algebraic varieties with a pencil of del Pezzo
surfaces, Izvestiya: Mathematics {\bf 62}:1 (1998), 115-155.
\vspace{0.3cm}

\noindent 8. Pukhlikov A.V., Birational automorphisms of Fano
hypersurfaces, Invent. Math. {\bf 134} (1998), no. 2, 401-426.
\vspace{0.3cm}

\noindent 9. Pukhlikov A.V., Birationally rigid Fano double
hypersurfaces, Sbornik: Mathematics {\bf 191} (2000), No. 6,
101-126. \vspace{0.3cm}

\noindent 10. Pukhlikov A.V., Birationally rigid Fano fibrations,
Izvestiya: Mathematics {\bf 64} (2000), 131-150. \vspace{0.3cm}

\noindent 11. Pukhlikov A.V., Certain examples of birationally
rigid varieties with a pencil of double quadrics. Journal of Math.
Sciences. 1999. V. 94, no. 1, 986-995.
 \vspace{0.3cm}

\noindent 12. Pukhlikov A.V., Birational automorphisms of
algebraic varieties with a pencil of double quadrics. Mathematical
Notes. {\bf 67} (2000), 241-249. \vspace{0.3cm}

\noindent 13. Pukhlikov A.V., Essentials of the method of maximal
singularities, in ``Explicit Birational Geometry of Threefolds'',
London Mathematical Society Lecture Note Series {\bf 281} (2000),
Cambridge University Press, 73-100. \vspace{0.3cm}

\noindent 14. Pukhlikov A.V., Birationally rigid Fano complete
intersections, Crelle J. f\" ur die reine und angew. Math. {\bf
541} (2001), 55-79. \vspace{0.3cm}

\noindent 15. Pukhlikov A.V., Birationally rigid iterated Fano
double covers. Izvestiya: Mathematics. {\bf 67} (2003), no. 3,
555-596. \vspace{0.3cm}

\noindent 16. Pukhlikov A.V., Birationally rigid varieties with a
pencil of Fano double covers. I. Sbornik: Mathematics {\bf 195}
(2004). Available as a preprint MPI-2003, no. 102.  \vspace{0.3cm}

\noindent 17. Sarkisov V.G., Birational automorphisms of conic
bundles, Izv. Akad. Nauk SSSR, Ser. Mat. {\bf 44} (1980), no. 4,
918-945 (English translation: Math. USSR Izv. {\bf 17} (1981),
177-202). \vspace{0.3cm}

\noindent 18. Sarkisov V.G., On conic bundle structures, Izv.
Akad. Nauk SSSR, Ser. Mat. {\bf 46} (1982), no. 2, 371-408
(English translation: Math. USSR Izv. {\bf 20} (1982), no. 2,
354-390). \vspace{0.3cm}

\noindent 19. Sobolev I. V., On a series of birationally rigid
varieties with a pencil of Fano hypersurfaces. Mat. Sb. {\bf 192}
(2001), no. 10, 123-130 (English translation in Sbornik: Math.
{\bf 192} (2001), no. 9-10, 1543-1551). \vspace{0.3cm}

\noindent 20. Sobolev I. V., Birational automorphisms of a class
of varieties fibered into cubic surfaces. Izv. Ross. Akad. Nauk
Ser. Mat. {\bf 66} (2002), no. 1, 203-224. \vspace{0.3cm}

\noindent 21. Fulton W., Intersection Theory, Springer-Verlag,
1984. }

\begin{flushleft}
{\it e-mail}: pukh@mi.ras.ru, pukh@liv.ac.uk
\end{flushleft}

\end{document}